\documentclass{article}
\usepackage{euscript,amsfonts,amssymb,amsmath,multicol,amsthm,slashbox,longtable}
\usepackage[all]{xy}
\usepackage{graphicx,epsfig}
\author{A. G. Gorinov}
\title{Division theorems for the rational cohomology of some discriminant complements}
\date{}

\newcommand{\Z}{{\mathbb Z}}
\newcommand{\C}{{\mathbb C}}
\newcommand{\R}{{\mathbb R}}
\newcommand{\Q}{{\mathbb Q}}
\newcommand{\dd}{\underline{d}}
\newcommand{\aaa}{\mathbf{a}}
\newcommand{\bb}{\mathbf{b}}
\newcommand{\cc}{\mathbf{c}}
\newcommand{\ee}{\mathbf{e}}
\newcommand{\m}{\mathbf{m}}
\newcommand{\W}{\mathbf{W}}
\newcommand{\X}{\mathbf{X}}
\newcommand{\Y}{\mathbf{Y}}

\newcommand{\LCM}{\mathrm{LCM}}
\newcommand{\x}{\mathbf x}
\newcommand{\y}{\mathbf y}
\newtheorem{theorem}{Theorem}
\newtheorem*{theoremx}{Theorem {\thevb}${}^\prime$}
\newtheorem{lemma}{Lemma}
\newtheorem{Prop}{Proposition}

\newtheorem{corollary}{Corollary}
\DeclareMathOperator\lk{lk}
\DeclareMathOperator\rk{rk}
\DeclareMathOperator\Gr{Gr}
\DeclareMathOperator\Sing{Sing}
\DeclareMathOperator\codim{codim}

\newcommand{\PGL}{\mathrm{PGL}}
\newcommand{\GL}{\mathrm{GL}}
\newcommand{\SO}{\mathrm{SO}}
\newcommand{\Mat}{\mathrm{Mat}}
\DeclareMathOperator\tot{tot}
\DeclareMathOperator\Aut{Aut}
\begin{document}
\maketitle
\begin{abstract}
The main purpose of this paper is to show that the mixed Hodge polynomial of the ``space of equations'' for smooth complete intersections of given multidegree
in $\mathbb{C} P^n$
is divisible by the mixed Hodge polynomial of the group $\mathrm{GL}_{n+1}(\mathbb{C})$,
the quotient being the mixed Hodge polynomial of the corresponding quotient space.
As a by-product of the method used in the proof, we obtain expressions divisible by the order the automorphism group of any smooth projective hypersurface of given dimension and degree.
\end{abstract}
\section{Introduction and main results}\label{intro}
\label{nkdd}Let $n$ and $k$ be integers satisfying $1\leq k\leq n+1$.
%and let $\kk$ be an algebraically closed field of characteristic zero.
Set $\dd=(d_1,\ldots,d_k)$ to be a collection
of integers such that $2\leq d_1\leq\cdots\leq d_k$. Denote by \label{PiSigma}$\Pi_{\dd,n}$ the $\C$-vector
space of all $k$-tuples $(f_1,\ldots ,f_k)$, where $f_i,i=1,\ldots, k$, is a homogeneous polynomial
in $n+1$ variables of degree $d_i$ with coefficients in $\C$. For every $(f_1\ldots ,f_k)\in\Pi_{\dd,n}$
denote by \label{SING}$\Sing (f_1,\ldots ,f_k)$ the projectivisation of the set of all $x\in\C^{n+1}\setminus\{ 0\}$ such that
\begin{itemize}
\item $f_i(x)=0,i=1,\ldots,k$,
\item the gradients of $f_i,i=1,\ldots,k$ at $x$ are linearly dependent.
\end{itemize}
Set $\Sigma_{\dd,n}$ to be the
subset of $\Pi_{\dd,n}$ consisting of all $(f_1,\ldots ,f_k)$ such that
$\Sing (f_1,\ldots ,f_k)\neq\varnothing$.
%For any $(f_1,\ldots,f_k)\in\Pi_{\dd,n}$ set $\langle f_1\ldots,f_k\rangle$ to be the ideal generated by $f_1,\ldots ,f_k$.
If $(f_1,\ldots ,f_k)\in\Pi_{\dd,n}\setminus\Sigma_{\dd,n}$, then the subvariety $X$ of $\C P^n$ defined by $f_1=\cdots=f_k=0$ is smooth, and $f_1,\ldots,f_k$ generate the homogeneous ideal of $X$. For this reason the space~$\Pi_{\dd,n}\setminus\Sigma_{\dd,n}$ can be viewed as the space of equations
for some smooth complete intersections of multidegree $\dd$ in $\C P^n$.
The case $k=n+1$ does not does not quite agree with this interpretation (it would correspond to ``empty complete intersections''), but we include it nonetheless, since we shall need it later as a starting point for some
computations.

%Let $\ell(\dd)$ be the number of different items in $\dd$. Let $(m_1,\ldots,m_{\ell(\dd)})$ be the sequence
%of integers such that for any $i=1,\ldots,\ell(\dd)$
%$$d_{1+\sum_{j<i}m_j}=\cdots=d_{m_i+\sum_{j<i}m_j},d_{\sum_{l\leq i}m_j}<d_{\sum_{j\leq i}m_j+1}\mbox{(for $i<\ell(\dd)$)}.$$
%%Denote by $\langle f_1,\ldots,f_l\rangle$ the ideal generated by polynomials $f_1,\ldots,f_l$.
%There exists an algebraic group $\mathbf{G}_{\dd,n}$ acting on $\Pi_{\dd,n}\setminus\Sigma_{\dd,n}$
%such that
%
%\begin{itemize}\label{ideals}
%\item $\langle f_1,\ldots ,f_k\rangle=\langle g_1,\ldots ,g_k\rangle$, iff $(f_1,\ldots ,f_k)$ and $(g_1,\ldots ,g_k)$
%belong to the same $\mathbf{G}_{\dd,n}$-orbit,
%\item a maximal compact subgroup of $\mathbf{G}_{\dd,n}$ is isomorphic to $U_{m_1}\times\cdots U_{m_{\ell(\dd)}}$.
%\end{itemize}
%(We shall describe this group and the action in the next section.)
%The group $\GL_{n+1}(\C)\times\mathbf{G}_{\dd,n}$ contains a subgroup isomorphic to $\C^*$ that acts identically on $\Pi_{\dd,n}\setminus\Sigma_{\dd,n}$ (see the next
%section); denote by ${\cal{G}}_{\dd,n}$ the corresponding quotient group.
%

The group $\GL_{n+1}(\C)$ acts on $\Pi_{\dd,n}$ in an obvious way:
\begin{equation}\label{firstaction}\GL_{n+1}(\C)\times \Pi_{\dd,n}\ni (A,(f_1,\ldots,f_k))\mapsto (f_1\circ A,\ldots, f_k\circ A);\end{equation}
this action preserves $\Sigma_{\dd,n}$ (and hence, $\Pi_{\dd,n}\setminus\Sigma_{\dd,n}$).
%We shall write $\Pi_{\dd,n}$ and
%$\Sigma_{\dd,n}$ instead of $\Pi_{\dd,n}(\C)$ and $\Sigma_{\dd,n}(\C)$ respectively.
%In the sequel we work almost entirely over $\C$, however some of the results can be extended to the case of abritrary $\kk$ completely
%for free, which is why we introduced the notation in a more general setting.

The main purpose of the paper is to prove the following theorem.

\begin{theorem}\label{main}
Suppose $\dd\neq (2)$. Then the geometric quotient of $\Pi_{\dd,n}\setminus\Sigma_{\dd,n}$ by $\GL_{n+1}(\C)$ exists, and
the Leray spectral sequence of the corresponding quotient map degenerates over $\Q$ (or modulo a sufficiently large prime)
at the second term.
\end{theorem}
This theorem generalises a recent result of J. Steenbrink and C. Peters for the case $k=1$ \cite{stepet}.
Our general strategy will be the same as in \cite{stepet}; however, the details will be different and more elementary (or so we hope).

For a complex algebraic variety $V$, we define the {\it mixed Hodge polynomial} of $V$ to be
$$P_{\mathrm{mHdg}}(V)=\sum_{n,p,q}t^nu^pv^q\dim_\C(\Gr^p_F\Gr_{p+q}^WH^n(V,\C)).$$
By setting in this expression $v=u$, respectively, $u=v=1$, we get the Poincar\'e-Serre polynomial, respectively, the Poincar\'e polynomial, of $V$.
The {\it mixed Hodge polynomial of $V$ with compact supports} (which we denote by $P_{\mathrm{mHdg},c}$) is obtained by replacing $H^n$ by $H^n_c$ in the definition of $P_{\mathrm{mHdg}}$;
by specialising $P_{\mathrm{mHdg},c}$ at $t=-1$ we get the Serre characteristic of~$V$.

An easy corollary of theorem \ref{main} is
\begin{corollary}
We have $$P_{\mathrm{mHdg}}(\Pi_{\dd,n}\setminus\Sigma_{\dd,n})=P_{\mathrm{mHdg}}(\GL_{n+1}(\C))\cdot P_{\mathrm{mHdg}}((\Pi_{\dd,n}\setminus\Sigma_{\dd,n})/
\GL_{n+1}(\C)),$$
$$P_{\mathrm{mHdg,c}}(\Pi_{\dd,n}\setminus\Sigma_{\dd,n})=P_{\mathrm{mHdg,c}}(\GL_{n+1}(\C))\cdot P_{\mathrm{mHdg,c}}((\Pi_{\dd,n}\setminus\Sigma_{\dd,n})/
\GL_{n+1}(\C)).$$
\end{corollary}
%The elements of the quotient $(\Pi_{\dd,n}\setminus\Sigma_{\dd,n})/\mathbf{G}_{\dd,n}$ parametrise smooth
%complete intersections in $\C P^n$ of multidegree $\dd$.
%The group $\PGL_{n+1}(\C)$ acts on $(\Pi_{\dd,n}\setminus\Sigma_{\dd,n})/\mathbf{G}_{\dd,n}$. An easy consequence of theorem \ref{main}
%is
%
%\begin{theorem}\label{main1}
%The geometric quotient of $(\Pi_{\dd,n}\setminus\Sigma_{\dd,n})/\mathbf{G}_{\dd,n}$ by $\PGL_{n+1}(\C)$ exists, and
%the rational Leray spectral sequence of the corresponding quotient map degenerates
%at the term $E_2$.
%\end{theorem}
%
By the Leray-Hirsch principle, in order to prove theorem \ref{main}, it suffices to construct global cohomology classes on $\Pi_{\dd,n}\setminus\Sigma_{\dd,n}$
(over $\Q$ or modulo a prime $p,p\gg 0$) such that their pullbacks under any orbit map generate the cohomology of the group $\GL_{n+1}(\C)$ (as a topological space).
We realise such classes as linking numbers with some
natural subvarieties of~$\Sigma_{\dd,n}$.

It turns out that in our situation working with integer coefficients is just a little bit more difficult than with the rationals. However, taking this little
extra effort pays off, since it enables one to determine explicitly which multiple of the generator of the highest cohomology group of $\GL_{n+1}(\C)$ comes from
the cohomology of $\Pi_{\dd,n}\setminus\Sigma_{\dd,n}$ via an orbit map. This (together with some simple computations) implies the following results:

\newcounter{vb}
\begin{theorem}\label{upperboundvector}
Let $d$ be an integer $>2$. Then the
order of the subgroup of $\GL_{n+1}(\C)$ consisting of the transformations that fix $f\in\Pi_{(d),n}\setminus\Sigma_{(d),n}$ divides
$$\prod_{i=0}^n((-1)^{n-i}+(d-1)^{n-i+1})(d-1)^i.$$
\end{theorem}
\setcounter{vb}{\value{theorem}}

Actually, we prove in section \ref{sec5} an analogous statement for arbitrary $\dd$ (theorem \thevb${}^{\prime}$), but the resulting formula is a bit messy
(and was therefore banned from the introduction).

\begin{theorem}\label{upboundprojective}
The order of the subgroup $\PGL_{n+1}(\C),n\geq 1,$ consisting of the transformations that preserve a smooth hypersurface of degree $d>2$ divides
\begin{equation}\label{bound}
\frac{1}{n+1}\prod_{i=0}^{n-1}\frac{1}{C^i_{n+1}}((-1)^{n-i}+(d-1)^{n-i+1})\LCM (C^i_{n+1}(d-1)^i,(n+1)(d-1)^n).
\end{equation}
%In particular, if $d\neq n+1,n\geq 3$, then the order of the automorphism group of $X$ divides~(\ref{bound}).
\end{theorem}
(Here $\LCM$ stands for the least common multiple.)

By the Lefschetz principle, the statements of theorems \ref{upperboundvector} and \ref{upboundprojective} are in fact true over
any algebraically closed field of characteristic 0.

In the case $n=1$, theorem \ref{upboundprojective} is equivalent to saying that the order of the group of linear fractional transformations that preserve
a (given) subset of $d>2$ points of $\C P^1$ divides $d(d-1)(d-2)$; this should be easy to prove directly. Notice that if $d,d-1$ and $d-2$ are pairwise coprime, one can not expect to have a
strictly stronger result.

For curves in $\C P^2$, surfaces in $\C P^3$ and threefolds in $\C P^4$ the expression (\ref{bound}) amounts to
\begin{equation}\label{krivye}d^2(d-1)^4(d^2-3d+3)(d-2),\end{equation}
$$\frac{1}{3}d^3(d-1)^8(d^3-4d^2+6d-4)(d^2-3d+3)(d-2)\LCM(3,2(d-1))$$
and
\begin{multline*}\frac{1}{4}d^4(d-1)^{13}(d^4-5d^3+10d^2-10d+5)(d^3-4d^2+6d-4)(d^2-3d+3)(d-2)\\
\cdot\LCM(2,(d-1)^2)\LCM(2,d-1)\end{multline*}
respectively.

In the case of plane curves,
a result similar to formula (\ref{krivye}) was obtained in~\cite{af} by P. Aluffi and C. Faber by studying the degrees of the $\PGL_3(\C)$-orbits of smooth curves.
Namely, in \cite{af} the authors produce, given a smooth plane curve $C$, a certain expression depending on the degree and the nature of flexes of $C$ and divisible by the order of the stabiliser $\subset\PGL_3(\C)$ of $C$.

If $n\geq 3,d\geq 3$ and $(d,n)\neq(4,3)$, then any automorphism of a smooth hypersurface of degree $d$ in $\C P^n$ is known \cite[theorem 2]{mm} to be the restriction of a projective
transformation, so in these
cases theorem \ref{upboundprojective} implies that the order of the full automorphism group divides~(\ref{bound}).
%It is well known (see, e.g., \cite[page ...]{gh})
%that any automorphism of a smooth hypersurface of degree $d\neq n+1$ in $\C P^n, n\geq 3,$ is the restriction of a projective
%transformation of $\C P^n$, hence theorem \ref{upboundprojective} implies that the order of the automorphim group of such a hypersurface
%divides (\ref{bound}).

The expression (\ref{bound}) is majorated by $d^{\frac{3}{2}n(n+1)}(n+1)^{n-1}$; since (\ref{bound}) is divisible
by the order of the projective automorphism group of any smooth hypersurface of degree $d$ in $\C P^n$, it can hardly be expected to be a sharp bound.
Indeed, smaller bounds are known;
the best one known to the author is
%\begin{equation}\label{hs}b(n)\begin{vmatrix}\sum_{i=1}^n(-1)^{i+1}d^iC_{n+1}^{n-i}\end{vmatrix}\end{equation}
\begin{equation}\label{hs}J(n+1)d^n\end{equation}
given by
A. Howard and A. J. Sommese \cite{hoso} (here $J$ is the Jordan function, i.e., $J(m)$ is the minimal integer such that any finite subgroup of $\GL_m(\C)$
contains a normal Abelian subgroup of index $\leq J(m)$;
B. Weisfeiler proved \cite{weis} that $J(m)\leq (m+1)!m^{a\ln m +b}$ for some $a,b\in\R$). However, theorem \ref{upboundprojective} gives additional information
on the orders of automorphism groups; in a sense,
asymptotically as $d\to\infty$, it provides much more restrictions than (\ref{hs}), since the number of divisors
of $x\in\Z$ grows more slowly than any power of $x$ as $x\to\infty$ (see, e.g., \cite[theorem 317]{hardy}).
More on theorem \ref{upboundprojective} can be found in section \ref{discussion}.

Upper bounds for automorphism groups of arbitrary varieties of general type are given in \cite{sz}. See also G. Xiao's papers \cite{xiao1,xiao2}
where it is shown that the automorphism group of a surface $S$ of general type contains at most $1764 (K_S)^2$ elements (this bound amounts to
$1764d(d-4)^2$ for surfaces of degree $d>4$ in $\C P^3$).
%Theorem \ref{upboundprojective} improves (when applicable, of course) the general upper bounds by E. Szab\'o
%\cite{sz}.
%For smooth projective complex hypersurfaces, general results stronger than those given by \cite[theorem 4]{sz} may be known. The best upper bound
%one can hope for is a polynomial of degree $n$ in $d$ (for any fixed $n$). However, I was unable to find such results in the literature.
%
%For curves and surfaces there exist the classical Hurwitz bound ($84(g-1)$) and G. Xiao's bound ($1764 K^2$, see \cite{xiao1,xiao2}) respectively (the latter amounts to
%$1764d(d-n-1)^2$ for surfaces of degree $d>4$ in $\C P^3$). However, even in these cases, theorem \ref{upboundprojective} gives additional information
%on the orders of automorphism groups; in a sence, asymptotically, it provides more restrictions than the bounds of Hurwitz and G. Xiao, since the number of divisors
%of $m$ grows more slowly than any power of $m$ as $m\to\infty$ (see, e.g., \cite[theorem 317]{hardy}).
%More on the comparison of theorem \ref{upboundprojective} to the previously known results can be found in section \ref{discussion}.
%
%certain multiples of some of these global classes descend into
%$(\Pi_{\dd,n}\setminus\Sigma_{\dd,n})/\mathbf{G}_{\dd,n}$, which gives immediately an upper bound on the ambient automorphism group of a smooth
%complete intersection of mutlidegree $\dd$ in $\C P^n$.
%

The idea of the proof of theorem \ref{main} presented here came from the following remark. The first columns of the Vassiliev spectral sequences
that compute the Borel-Moore homology of the determinant varieties \cite{vas1} (i.e., the spaces of degenerate $(n+1)\times (n+1)$-matrices)
and of the discriminant varieties $\Sigma_{\dd,n}$ (see \cite{vas2,quintics,tommasi}) coincide up to a dimension shift. This paper may be
viewed as an attempt to understand the relationship between the corresponding cohomology classes.

The (rest of the) paper is organised in the following way. In section \ref{sec2} we introduce some
notation and formulate and/or prove several preliminary results. Then we study in sections \ref{sec3}
and~\ref{sec4} the way homology classes of the $\GL_{n+1}(\C)$-orbits in $\Pi_{\dd,n}\setminus\Sigma_{\dd,n}$ are linked with certain
subvarieties of $\Sigma_{\dd,n}$; the results of these sections are then used in section \ref{sec5}
to prove theorems~\ref{main}-\ref{upboundprojective}.
In section \ref{discussion} we list the values of (\ref{bound}) for small $n$ and $d$ and discuss several particular cases of theorem \ref{upboundprojective}. We also
give there an analogue of theorem \ref{upboundprojective} for the groups of deck transformations of ramified coverings $\C P^n\to\C P^n$ and an application of theorem~\ref{main} to real algebraic geometry and discuss some open questions. In the end we give an index of
some non-self-explanatory notation used throughout the paper.

I am grateful to Joseph Steenbrink for the interest, encouragement and many useful conversations; I would also like to thank Tatiana Bandman,
Pierre Vogel and Mikhail Zaidenberg
for useful discussions and correspondence.
\section{Notation and preliminaries}\label{sec2}
In the sequel the coefficients of all (co)homology groups are assumed to be integer, unless stated otherwise. Whenever different sign choices are possible (e.g., in the definition of the $\frown$-product), we shall use the classical conventions (e.g., as in \cite{fuchs}).
%
%In what follows, $\C^m$ will be considered as a vector space and its elements as column vectors.
%
%If $V\subset\C^m$ is an affine variety, we denote the ideal of $V$ by $I(V)$.

Any topological space considered in \ref{BM1} and \ref{BM2} is assumed to have the form ``a finite $CW$-complex minus a point''. Notice that any
complex algebraic variety satisfies this condition.

All algebraic varieties that we shall consider will be
defined over $\C$, unless stated otherwise. As usual, we denote by $I(V)$ the ideal of an affine variety
$V$.

We shall denote by \label{tot}$\tot(\xi)$ the total space of a vector bundle $\xi$.

The first Chern class of the cotautological bundle on $\C P^m$ will be called the {\it canonical generator} of $H^*(\C P^m)$.

In the sequel we consider elements of $\C^m$ as column vectors, unless stated otherwise.

When talking about the cohomology of an algebraic variety or a Lie
group we shall always mean the cohomology of
the corresponding topological space.

\subsection{Borel-Moore homology and linking numbers}\label{BM1}
Let $X$ be a topological space. We denote by $\bar X$ the one-point compactification of $X$.
The Borel-Moore
homology groups of $X$ will be denoted by $\bar H_*(X)$. These groups can be viewed
either as the homology groups
of the complex of locally finite singular chains or as the homology groups of $\bar X$ modulo the infinity.
%We shall say that
%$X$ {\it orientable}, if $\bar X$ is homeomorphic to a finite $CW$-complex of dimension $m$, and $\bar H_m(X)\cong\Z$, in which case we shall call $X$
%{\it oriented}, if one of the two isomorphisms $\Z\to \bar H_m(X)$ is chosen; the image of 1 under this preferred isomorphism
%will be called the {\it fundamental class} of $X$ and denoted by $[X]$.

We shall use the symbol $H_{/\mathrm{tor}}^*$ to denote the integer cohomology groups modulo torsion.

Let $M$ be a smooth oriented manifold of dimension $p$, let $X\subset M$ a closed subspace, and let $c\in\ker (\bar H_{p-q}(X)\to\bar H_{p-q}(M))$.
Suppose $H_{q-1}(M)=0$. Then the group $H^{q-1}(M)\cong\bar H_{p-q+1}(M)$ is finite, and
$c$ defines a unique element of $\bar H_{p-q+1}(M\setminus X)/\mbox{torsion}\cong H_{/\mathrm{tor}}^{q-1}(M\setminus X)$. This element will be
called {\it the linking number with $c$ in $M$} and denoted by
\label{llk}$\lk_{c,X,M}$.
Notice that if, in addition to $H_{q-1}(M)=0$ we also have $H^{q-1}(M)=0$, which will be the case in all applications that we have 
in mind, then $\lk_{c,X,M}$ can be naturally viewed as a cohomology class $\in H^{q-1}(M\setminus X)$.

Here is another equivalent definition: suppose that $c$ is represented by a smooth singular chain $\tilde c$, and consider the function
$H_{q-1}(M\setminus X)\to\Z$ defined
as follows: take a cycle $\in H_{q-1}(M\setminus X)=\ker(H_{q-1}(M\setminus X)\to H_{q-1}(M))$, represent it by
a smooth singular chain $z$, find a smooth singular chain $w$ in $M$ that is bounded by $z$ and transversal to $\tilde c$, and calculate
the intersection index $\#(w, \tilde c)$. This function also defines a unique element
of $H_{/\mathrm{tor}}^{q-1}(M\setminus X)$, which coincides with $(-1)^q\lk_{c,X,M}$.

If $Y$ is a closed supspace of $M$ containing $X$, then an easy check shows that $\lk_{c,X,M}$ restricted to $M\setminus Y$ is $\lk_{c',Y,M}$, where
$c'$ is the image of $c$ in $\bar H_{p-q}(Y)$:
\begin{equation}\label{diagxy}\xymatrix{H^{q-1}(M\setminus Y) & \bar H_{p-q+1}(M\setminus Y)\ar[l]_-{\mbox{\scriptsize  Poincar\'e}}\ar[r]^-\cong &
\bar H_{p-q+1}(M,Y)\ar[r] & \bar H_{p-q}(Y)\\
H^{q-1}(M\setminus X)\ar[u] & \bar H_{p-q+1}(M\setminus X)\ar[l]_-{\mbox{\scriptsize Poincar\'e}}\ar[r]^-\cong\ar[u] &
\bar H_{p-q+1}(M,X)\ar[r]\ar[u] & \bar H_{p-q}(X)\ar[u]}\end{equation}
In particular, if $V_1\subset V_2\subset \C^m$ are irreducible
affine subvarieties, then the restriction of $\lk_{[V_1],V_1,\C^m}$ to
$\C^m\setminus V_2$ is Alexander dual to the image of $[V_1]$ in $\bar H_*(V_2)$.

In the sequel, if $X\subset Y\subset M$, we shall often write $\lk_{c,X,M}$ instead of $\lk_{c,X,M}|_{M\setminus Y}$.
We shall also drop $X$ from the notation whenever it is
clear which space $X\subset M$ we are considering.

Suppose that $X_1$ and $X_2$ are topological spaces, $f:X_1\to X_2$ is a locally trivial fibration with fibre $F$, which we assume to be a smooth orientable manifold. Suppose that $f$ is homologically trivial, and introduce an orientation on $F$. There exists a {\it preimage map} $f^+_*:\bar H_*(X_2)\to\bar H_{*+\dim_\R F}(X_1)$ defined as follows: we associate
to an element $a\in \bar H_*(X_2)$ the image in $\bar H_*(X_1)$ of
what has remained in $\bar E^\infty$ of the class $a\otimes[F]\in\bar E^2$ (where $(\bar E^r)$ is the Leray
spectral sequence converging to the Borel-Moore homology of $X_1$).

If $X_2$ is a smooth manifold, and $a$ is represented by a smooth singular chain $\tilde a$, then $f_*^+(a)$ can be represented by the (appropriately triangulated) preimage of the support of~$\tilde a$ (hence the term ``preimage map'').
%In particular, if $X_2$ is a smooth manifold, another (equivalent) construction of $f^+_*$ can be given:
%take an element $a\in H_*(X_2)$, represent is by a (real) possibly singular oriented submanifold $X$ with singularities of (real) codimension $\geq 2$, %and
%set $f^+_*(a)=[p^{-1}(X)]$ (hence the term ``preimage map'').

The basic properties of linking numbers that we shall need can be summarised as follows.

\begin{Prop}
\begin{enumerate}\label{lk}\item
Suppose that $M_1$ and $M_2$ are smooth algebraic varieties, $f:M_1\to M_2$ is a locally trivial holomorphic fibration with fibre $F$, and let $N_2$ be a (closed) subvariety of $M_2$. Set $N_1=f^{-1}(N_2)$. Assume that $f$ is homologically trivial. We have then $f^*(\lk_{c,N_2,M_2})=\lk_{f^+_*(c),N_1,M_1}$ for any $c\in\bar H_*(N_2)$ whose image in $\bar H_*(M_2)$ is zero.
%$N^1_2,\ldots,N^l_2\subset M_2$ a collection of
%?closed irreducible pairwise distinct subvarieties
%of codimension $q$. Set $N^i_1=f^{-1}(N^i_2),i=1,\ldots ,l$. Suppose that $H_{2q-1}(M_1)=H_{2q-1}(M_2)=0$ and that
%the image of $\sum_{i=1}^l\alpha_i[N^i_j]$ in $\bar H_*(M_j),j=1,2$ is zero, where $\alpha_1,\ldots,\alpha_l$ are integers
%(notice that the image of any $[N^i_j]$ in $\bar H_*(M_j),j=1,2$ might well be nonzero).
%Then the following holds in
%$H^*_{/\mathrm{tor}}(M_1\setminus\bigcup N^i_1)$: $$f^*(\lk_{\sum_{i=1}^l\alpha_i[N^i_2],M_2})=\lk_{\sum_{i=1}^l\alpha_i[N^i_1],M_1}.$$
%\item Suppose that $X\subset\R^m$ is a closed oriented $CW$-complex. Then $\lk_{X,\R^m}$ is Alexander dual to $[X]$.
\item Suppose that $V\subset\C^m$ is an irreducible  affine subvariety, and let $E\subset\C^m$ be an affine plane
such that $\codim_{\C^m}V=\codim_EV\cap E$. Let $C_1,\ldots,C_l$ be the components of $V\cap E$ of maximal dimension. Denote by
$\alpha_i$ the intersection multiplicity of $V$ and $C$ along $C_i,i=1,\ldots,l$. Then in $H^*(E\setminus V)$ the following holds: $$\lk_{[V],\C^m}|_{E\setminus V}=\sum_{i=1}^l\alpha_i \lk_{[C_i],E}.$$
\item Suppose that $V\subset\C^m$ is an irreducible affine variety, set $l=\codim_{\C^m}V$, and let $F:\C^l\to\C^m$ be a polynomial mapping
whose restriction to the unit ball $U\subset\C^l$ is an embedding and such that $F^{-1}(V)\cap U=\{0\}$.
Then the pullback under $F|_U$ of $\lk_{[V],\C^m}$ is the canonical generator of $H^{2l-1}(U\setminus\{0\})$ times
intersection multiplicity $\mu$ of $F(U)$ and $V$ at $F(0)$. We have
$$\mu=\dim_\C{\cal O}_{0,\C^l}/F^{*}(I(V))_{\mathfrak{m}_0},$$ where $\mathfrak m_0$ is the ideal formed
by the polynomials in ${\cal O}_{\C^l}$ that vanish at the origin. If moreover $F^{-1}(V)=\{0\}$, then
$$\mu=\dim_\C({\cal O}_{0,\C^l}/F^{*}(I(V))_{\mathfrak{m}_0})=\dim_\C({\cal O}_{\C^l}/F^{*}(I(V))).$$
%(Notice that if the variety $V$ is conical, i.e., if $x\in V$ implies $\lambda x\in V$ for any $\lambda\in\C$,
%then $I$ is just the ideal of $V$.)
\end{enumerate}
\end{Prop}
$\clubsuit$

\begin{corollary}\label{linear}
If $f:E_1\to E_2$ is a surjective linear map of complex vector spaces, and $V\subset E_2$ is an irreducible affine subvariety, then
$$f^*(\lk_{[V],E_2})=\lk_{[f^{-1}(V)],E_1}.$$
\end{corollary}
$\clubsuit$

\subsection{Some facts about vector bundles}\label{BM2}
All propositions in this subsection are standard exercises, but proofs are given nonetheless for the sake of completeness.

\subsubsection{The image of the zero section in the Borel-Moore homology}

\begin{Prop}\label{zerosec}
Let $\eta$ be a real oriented vector bundle of rank $l$ over a real smooth oriented manifold $X$
of dimension $d$, and let $Y\subset X$ be a oriented submanifold Poincar\'e dual to the Euler class
$e(\xi)$ of $\xi$. Set $E=\tot(\eta),E'=\tot(\eta|_Y)$. The image of $[X]$ in $\bar H_*([E])$ under the zero section embedding is equal to $[E']$.
\end{Prop}
(We orient $E$ and $E'$ by the usual rule ``first the fibre, then the base''.)
%This is a standard exercise, but we give a proof, just in case.

{\bf Proof.}
Equip $\eta$ with a Riemannian metric, and set $E_0$ and $E_0'$ to be the union of all elements of $E$, respectively, of $E'$, of length $\geq 1$.
Denote by $u\in H^l(E,E_0)$ the Thom class of $\eta$, and set $e'$ to be the restriction of $u$ to $E$. By definition, $e(\eta)$
is the image of $e'$ under the isomorphism
$H^l(E)\cong H^l(X).$

Let $\bar B$ be some contractible compact neighbourhood of $\infty$ in $\bar X$, and set $B=\bar B\setminus\{\infty\},Z=\tot
(\eta|_B)$. Due to the functoriality of the $\frown$-product, the following diagram is commutative.
%triples (E,\emptyset,Z) and (E,E_0,Z)
$$\xymatrix{
\bar H_{d-l}(Y)\ar[r]&\bar H_{d-l}(X)\ar[r]^\cong & H_{d-l}(X,B)\ar[r]^\cong & H_{d-l}(E,Z)\\
&\bar H_d(X)\ar[r]^-{\cong}\ar[d]&H_d(X,B)\ar[u]_-{\bullet\frown e(\eta)}\ar[r]^-{\cong}&H_d(E,Z)\ar[u]_-{\bullet\frown e'}\ar[d]\\
\bar H_d(E')\ar[r]&\bar H_d(E)\ar[rr]^-{\cong}&&H_d(E,E_0\cup Z)\ar@/_4pc/[uu]_-{\bullet\frown u}
}$$
The manifold $Y$ is chosen so that the images of $[X]$ and $[Y]$ in $H_{d-l}(E,Z)$ coincide. The Thom isomorphism
$\bullet\frown u:H_d(E,E_0\cup Z)\to H_{d-l}(E,Z)$ takes the image of $[E']$ to the image of $[Y]$, which proves the proposition.$\clubsuit$

\subsubsection{Degree of some varieties swept by linear subspaces}
\begin{Prop}\label{degreeprop}
Let $\eta$ be an holomorphic vector subbundle of rank $l$ of the trivial bundle $\C^N\times X$ over an irreducible projective variety $X$ of dimension $d$. Set
$V\subset\C^N$ be the union of all fibres of $\eta$, and denote by $A$ and $v$ the
matrix
$$\begin{pmatrix}
0&0&\cdots&0&-c_d(\eta)\\
1&0&\cdots&0&-c_{d-1}(\eta)\\
\hdotsfor{5}\\
0&\cdots&1&0&-c_2(\eta)\\
0&\cdots&0&1&-c_1(\eta)
\end{pmatrix}$$
and the vector $(0,\ldots,0,1)^T$ respectively. Let $w\in H^{2d}(X)$ be the last coordinate of $A^dv$. Then $w([X])\geq 0$; if $w([X])> 0$,
then $V$ has the expected dimension $d+l$, and $\deg V=w([X])$.
\end{Prop} 

{\bf Proof.} Set $Y$ to be the total space of the projectivisation of $\eta$. We shall consider $Y$ as a subset of $\C P^{N-1}\times X$.
%There exists an obvious inclusion $\imath:Y\to \C P^{N-1}\times X$.
Let $a$ be the canonical generator of $H^*(\C P^{N-1})$.
The proposition would follow if we manage to show that $a^{l+d-1}\otimes 1=(a^{l-1}\otimes 1) p^*(w)$, where
$p:Y\to X$ is the bundle projection.

Identify $H^*(X)$ with its image under $p^*$; the ring $H^*(Y)$ is generated over $H^*(X)$ by $b=a\otimes 1$ with the relation
$$b^l+\sum_{i=l-d}^{l-1}b^i c_{l-i}(\eta)=0.$$
Using this relation we obtain
$$b^{l+j-1}=\sum_{i=l-d}^{l-1}u_{i,j}b^i,$$
where $u_{i,j}\in H^{2(l+j-i-1)}(X)$ are such that $(u_{l-d,j},\ldots,u_{l-1,j})^T=A^jv$. In particular, $b^{l+d-1}=b^{l-1} w,$
since $u_{<l-1,d}=0$ for dimension reasons.
The proposition is proven.$\clubsuit$

\subsubsection{Chern classes of some bundles}
Let $d$ be a positive integer, and set $\eta$
%and $\eta_i,i=0,\ldots,n$,
to be the vector bundle over $\C P^n$ with total space
$$\{(f,x)\in\Pi_{(d),n}\times\C P^n\mid \Sing f\ni x\}$$
%and
%$$\{(f,x)\in\Pi_{(d),n}\times\C P^n\mid \frac{\partial f}{\partial x_i}(x')=0\mbox{ for any preimage $x'\in\C^{n+1}\setminus\{0\}$ of $x$}\}$$
%respectively. Denote by $\xi$ the constant bundle $\Pi_{(d),n}\times\C P^n$.
%Let us compute Chern classes of $\eta$.
%
%The quotient bundle $\xi/\eta$ is isomorphic to the direct sum
%$\bigoplus_{i=0}^n \xi/\eta_i.$ Let $l\subset\C P^n$ be a line.
%By taking appropriate constant sections of $\xi|_l$ one can construct sections of $(\xi/\eta_i)|_l,i=0,\ldots,n,$ with $d-1$ simple zeroes, hence
%$c_1(\xi/\eta_i)=(d-1)a$, where $a$ is the canonical generator of $H^*(\C P^n)$. This implies the following proposition.

\begin{Prop}\label{chern}
The total Chern class of $\eta$ is equal to $(1+(d-1)a)^{-n-1}$, where $a$ is the canonical generator of $H^*(\C P^n)$.
\end{Prop}
$\clubsuit$

\subsection{Tautological principal bundles}
We denote by $G_m(\C^N)$ the Grassmann manifold consisting of all $m$-dimensional complex vector subspaces of $\C^N$. Let $F_m(\C^N)$ be the total space of
the corresponding
tautological principal bundle, i.e.,
$$F_m(\C^N)=\{(E,(v_1,\ldots,v_m))\in G_m(\C^N)\times(\C^N)^{\times m}\mid \mbox{$v_1,\ldots ,v_m$ span $E$}\}.$$
Let $m_1,\ldots,m_l$ and $N_1,\ldots N_l$ be sequences of positive integers, and let $S$ be an irreducible subvariety of
$\prod_{i=1}^lG_{m_i}(\C^{N_i})$.
%Let $c$ be a collection of subvarieties of $\prod_{i=1}^mG_i(\C^{N_i})$ transversal to
%$S$ that have the form
%$$G_{m_1}(\C^{N_1})\times\cdots\times\mbox{a variety dual to some Chern class}\times\cdots\times G_{m_l}(\C^{N_l}).$$
Denote by $p$ the projection \begin{equation}
\label{taut}\prod_{i=1}^lF_{m_i}(\C^{N_i})\to\prod_{i=1}^lG_{m_i}(\C^{N_i}).\end{equation}
Set $G=\prod_{i=1}^l\GL_{m_i}(\C)$.
%There exists an obvious ``preimage'' map
%$p':H_*(S)\to\bar H_{*+s}(p^{-1}(S))$, where $s=\dim_\R G$. Namely, we associate
%to an element $a\in H_*(S)$ the image in $\bar H_*(F)$ of
%what has remained in $\bar E^\infty$ of the class $a\otimes[G]\in\bar E^2$ (here $(\bar E^r)$ is the Leray
%spectral sequence converging to the Borel-Moore homology of $F$). If $S$ is smooth, another (equivalent) construction can be given:
%take an element $a\in H_*(S)$, represent is by a (real) possibly singular submanifold $X$ with singularities of (real) codimension $\geq 2$, and
%set $p'(a)=[p^{-1}(X)]$.
%$p':H_*(\prod_{i=1}^mG_i(\C^{N_i}))\to H_{*+r}(\prod_{i=1}^mF_i(\C^{N_i}))$, where $s=\dim_\R(G)$.

\begin{Prop}\label{zzz}
For any $c\in H^{>0}(\prod_{i=1}^mG_{m_i}(\C^{N_i}))$, we have $p_*^+([S]\frown c|_S)=0$ (where $p^+_*$ is the preimage map defined in \ref{BM1}).
\end{Prop}
%If the variety $p^{-1}(S\cap V_1\cap\cdots\cap V_q)$ is reducible, its fundamental class is defined to be the sum of the fundamental classes of the components
%of maximal dimension.
In the sequel we shall need only a particular case of this proposition. However, we give a general version, since it could be useful in further applications
(and the proof takes two lines anyway).

{\bf Proof.} It is easy to check that $p^+_*([S]\frown c|_S)=p^+_*([S])\frown p^*(c)$; this class vanishes, since $p^*(c)=0$.$\clubsuit$

\subsection{The cohomology of $\GL_m(\C)$}
\label{ecm}Denote by $\ee_m$ the canonical generator of $H^{2m-1}(\C^{m}\setminus\{ 0\})$.
We shall call the map $\GL_i(\C)\to\GL_{j}(\C),i<j,$ given by $$A\mapsto\left(\begin{array}{cc}A&0\\0&1\end{array}\right)$$ the {\it canonical inclusion}.

Let $\cc^m_1,\ldots ,\cc^m_m,\cc_i^m\in H^{2i-1}(\GL_m(\C))$ be the system of multiplicative generators
such that
\begin{itemize}\item $\cc^m_m$ is the pullback of $\ee_m$ under
the map $A\mapsto\mbox{the last column of $A$}$,
\item the pullback of $\cc^m_i$ under the canonical inclusion
$\GL_{m-1}(\C)\subset\GL_m(\C)$ is $\cc_i^{m-1}$ for $m\geq 2,1\leq i<m$.
\end{itemize}
\label{ocm}We shall write $o(\cc^m_i)$ to denote a linear combination of monomials in $\cc^m_j,j<i$.

\begin{Prop}\label{glm}
\begin{enumerate}
\item The cohomology map induced by $\GL_m(\C)\ni A\mapsto A^T\in\GL_m(\C)$ takes $\cc^m_i$ to $(-1)^{i+1}\cc^m_i$.
\item The cohomology map induced by $\GL_m(\C)\ni A\mapsto A^{-1}\in\GL_m(\C)$ takes $\cc^m_i$ to $-\cc^m_i$.
\end{enumerate}
\end{Prop}

{\bf Proof:} an easy induction on $m$. The induction step is performed in either case as follows.
We know the image of $\cc^m_i,i<m$, since the restriction $H^*(\GL_m(\C))\to H^*(GL_{m-1}(\C))$ is injective in dimensions $<2m-1$.
Write the image of $\cc^m_m$ as $a\cc^m_m+o(\cc^m_m),a\in\Z$. Since $\cc_m^m$ generates (over $\Z$) the degree $2m-1$ part of
$\ker (H^*(GL_m(\C))\to\GL_{m-1}(\C))$, we see that $o(\cc^m_m)=0$. The coefficient $a$ is fixed by looking at the action on the
highest cohomology group; this action is
the identity, if the corresponding map restricted to the standard $\mathrm{U}_m\subset\GL_m(\C)$ is orientation-preserving, and minus identity otherwise.
$\clubsuit$

\subsection{Spaces of matrices and polynomials}
%
%We shall denote by $I(V)$ the ideal of an affine variety $V\subset\C^m$.

\subsubsection{Subvarieties of spaces of matrices}
We denote by \label{mat}$\Mat_{i,j}(\C)$ the space of all complex matrices with $i$ rows and $j$ columns.

Recall that $k\leq n+1$. Set \label{Wwkn}$$W_{k,n}=\{A\in\Mat_{n+1,k}(\C)\mid\rk A<k\},$$ 
%to be the variety formed by
%the matrices of rank $<k$, and denote by $\W_{k,n}$ the following
%subvariety of $\Mat_{n+1,k}(\C)\times\C^{n+1}$:
\label{Wkn}$$\W_{k,n}=\{(A,x)\mid \rk A<k,x^T A=0\}.$$

The codimension of $\W_{k,n}$ in $\Mat_{n+1,k}(\C)\times\C^{n+1}$ is $n+1$.

Set \label{XY}$$X_{k,n}=\{A\in\Mat_{n+1,k}(\C)\mid\rk A<k,\mbox{the last row of $A$ is zero}\}$$ and (for $k>1$)
\label{XXY}$$Y_{k,n}=\{A\in\Mat_{n+1,k}(\C)\mid\mbox{the last $k-1$ columns of $A$ form a matrix from $X_{k-1,n}$}\}.$$

Finally, set \label{XXYY}$\X_{k,n}=\{(A,x)\in\W_{k,n}\mid A\in X_{k,n}\}$ and \label{XXXYY}$\Y_{k,n}=\{(A,x)\in\W_{k,n}\mid A\in Y_{k,n}\}$.

\begin{Prop}\label{matrixlemma} Suppose that $n+1>k>1$.
Then the intersection of $\W_{k,n}$ and the vector subspace $\mathbf{E}$ of $\Mat_{n+1,k}(\C)\times\C^{n+1}$ formed by all $(A,x),
A=(a_{i,j})_{\genfrac{}{}{0pt}{}{0\leq i\leq n}{1\leq j\leq k}}$ such that
\begin{equation}\label{matrix}a_{n,2}=\cdots=a_{n,k}=0\end{equation}
is $\X_{k,n}\cup\Y_{k,n}$; the intersection multiplicity along each one of these components is 1.
\end{Prop}

{\bf Proof.} Set $E$ to be the vector subspace of $\Mat_{n+1,k}(\C)$ defined by (\ref{matrix}).
The variety $\W_{k,n}$ (respectively, $\X_{k,n}$ and $\Y_{k,n}$) contains an open dense subset that is (the total space of)
a vector bundle over the subset of $W_{k,n}$ (respectively, of $X_{k,n}$ and $Y_{k,n}$) formed by matrices of rank $k-1$. Hence, to prove the proposition, it is sufficient to show that the intersection multiplicity of $W_{k,n}$ and $E$ along both
$X_{k,n}$ and $Y_{k,n}$ is 1.

%The codimension of $X_{k,n}$ in $E$ is $n-k+2$.
Let $T_1$ be an affine mapping of the unit ball $U\subset\C^{n-k+2}$ to $E$ such that $T_1(U)$
intersects $X_{k,n}$ transversally at one smooth point; let us assume that this point is $T_1(0)$ and that
the last $k-1$ columns of $T_1(0)$ are linearly independent.

It is well known (see, e.g., \cite[theorem 2.10]{brvet}) that the ideal of $W_{k,n}$ in $\Mat_{n+1,k}(\C)$ is generated by all $k\times k$-minors; analogously, the ideal of $X_{k,n}$ in $E$ is generated by the polynomial
$A\mapsto a_{n,1}$ (where $A=(a_{i,j})_{\genfrac{}{}{0pt}{}{0\leq i\leq n}{1\leq j\leq k}}$) and the $k\times k$-minors involving the first $n$ rows. An immediate check shows that the localisations at 0 of the pullbacks of these ideals under $T_1$ coincide, and hence, the intersection multiplicity of $W_{k,n}$ and $E$ along $X_{k,n}$ is 1.

%Introduce coordinates $t_1,\ldots t_{n-k+2}$ in $U$.
%
%We have to show that the pullback under $T$ of the ideal $I$ generated by all $k\times k$-minors
%contains all variables $t_1,\ldots,t_{n-k+2}$. Write $T(t),t\in U$ as
%$({T}_{i,j})_{\genfrac{}{}{0pt}{}{1\leq i\leq n+1}{1\leq j\leq k}}$. We can assume that
%${T}_{n+1,1}=t_{n-k+2}$ and ${T}_{i,j},i<n+1$ depend only on $t_1,\ldots,t_{n-k+1}$. The image of $T|_{t_{n-k+2}=0}$ is transversal to $X_{k,n}$, hence
%the pullback $T^{-1}(I)$ contains $t_1,\ldots,t_{n-k+1}$. On the other hand, since the last $k-1$ columns of $T(0)$ are linearly independent,
%$T^{-1}(I)\ni pt_{n-k+2}$, where $p\in\C[t_1,\ldots,t_{n-k+1}]$ is a polynomial with nonzero constant term. This implies that $T^{-1}(I)$ contains %$t_{n-k+2}$
%as well.

The case of $Y_{k,n}$ can be considered in an analogous way. Namely, let $T_2:U\to E$ be an affine mapping such that
$T_2(U)$ intersects $Y_{k,n}$ transversally at one smooth point, which is $T_2(0)$. Assume that the left bottom item of $T_2(0)$ is nonzero. The ideal of $Y_{k,n}$ in $E$ is generated by all $(k-1)\times (k-1)$-minors involving the first $n$ rows and the last $k-1$ columns, and we proceed as above to conclude that the intersection multiplicity of $W_{k,n}$ and $E$ along $Y_{k,n}$ is also~1.$\clubsuit$

An alternative, albeit longer proof of proposition \ref{matrixlemma} can be given as follows. One could start by showing
%(by using the same arguments as in the proof of \cite[theorem 2.10]{brvet})
that the ideal of $\W_{k,n}$ in $\Mat_{n+1,k}(\C)\times\C^{n+1}$ is generated by the polynomials $$(A,x)\mapsto\mbox{ a $k\times k$-minor of $A$}$$ (where $A\in\Mat_{n+1,k}(\C),x\in\C^{n+1}$) and $n+1$ polynomials obtained from the relation $x^TA=0$. Generators for the ideals of $\X_{k,n}$ and $\Y_{k,n}$ in $\mathbf{E}$ can be found in an analogous way.  The intesection multiplicities of $\W_{k,n}$ and $\mathbf{E}$ along $\X_{k,n}$ and $\Y_{k,n}$ can then be computed directly.

%Let $T''$ be the composition of $T'$ and
%the embedding $\Mat_{n,k-1}(\C)\to E\subset\Mat_{n+1,k}(\C)$ given by $$A\mapsto\left(\begin{array}{cc}0&A\\1&0\end{array}\right).$$ We have %${T''}^{-1}(I)={T'}^{-1}(I)$;
%the latter ideal obviously contains any $t_i,i=1,\ldots, n-k+2$. The proposition is proven.$\clubsuit$
%%\subsubsection{The groups $\mathbf{G}_{\dd,n}$}
%Notice that there are the following obvious transformations of $\Pi_{\dd,n}$ that preserve the ideal generated by any $(f_1,\ldots,f_k)$:
%$$(f_1,\ldots,f_k)\mapsto (f_1,\ldots,af_i+gf_j,\ldots,f_k)$$
%(here $i\neq j$ are indices such that $d_i\geq d_j$, $a\in\C^*$, and $g$ is a homogeneous polynomial of degree $d_i-d_j$). Set $\mathbf{G}_{\dd,n}$ to be the
%group generated by all such transformations.
%
%\begin{Prop} The group $\mathbf{G}_{\dd,n}$ satisfies the conditions on page \pageref{ideals} and acts freely on $\Pi_{\dd,n}\setminus\Sigma_{\dd,n}$.\end{Prop}$\clubsuit$
%
%Notice that $\mathbf{G}_{\dd,n}$ is the semi-direct product of a normal Abelian translation group
%and the group $\GL_{m_1}(\C)\times\cdots\GL_{m_{\ell(\dd)}}(\C)$ acting as
%$$A\cdot(f_1,\ldots,f_k)\mapsto A(f_1,\ldots,f_k)^T$$
%(here $A$ is a $k\times k$-matrix made of diagonal blocks belonging to $\GL_{m_i}(\C),i=1,\cdots,\ell(\dd)$).
%
\subsubsection{Some subvarieties of $\Sigma_{\dd,n}$}
For any $X\subset\C P^n$ denote by \label{VddnX}$V_{\dd,n,X}$ the
subset of $\Pi_{\dd,n}$ consisting of all $(f_1,\ldots ,f_k)$ such that
$\Sing (f_1,\ldots ,f_k)\cap X\neq\varnothing$. Notice that $\Sigma_{\dd,n}=V_{\dd,n,\C P^n}$.

\begin{Prop}
For any $i=0,\ldots,n$, and any projective subspace $L$ of dimension $i$, $V_{\dd,n,L}$ is an irreducible affine subvariety of $\Pi_{\dd,n}$ of codimension $n-i+1$.
\end{Prop}
$\clubsuit$

In the sequel, whenever it is clear (or irrelevant), which projective subspace $L\cong\C P^i$ we are considering, we shall
write $V_{\dd,n,\C P^i}$ instead of $V_{\dd,n,L}$.

Set \label{a}$$\aaa_i^{\dd,n}=\lk_{[V_{\dd,n,\C P^{n-i+1}}],\Pi_{\dd,n}}.$$

%where $L$ is an projective subspace $\cong\C P^{n-i+1}$. Obviously, the class $\aaa_i^{\dd,n}$ does not depend on the choice of $L$.
%
%The group $\GL_{n+1}(\C)$ acts on the right on $\Pi_{\dd,n}$ by the formula
%$$\Pi_{\dd,n}\times\GL_{n+1}(\C)\ni(f_1,\ldots,f_k,A)\mapsto (f_1\circ A,\ldots, f_k\circ A).$$
If $f$ is a homogeneous polynomial in $n+1$ variables, $x\in\C^{n+1}$, we set
$df|_x$ to be the vector $$\left(\frac{\partial f}{\partial x_0}(x),\ldots,\frac{\partial f}{\partial x_n}(x)\right)^T;$$
notice that if $A\in\GL_{n+1}(\C)$, and $g=f\circ A$, then
\begin{equation}\label{gradient}
dg|_x=A^Tdf|_{Ax}.
\end{equation}

%The action of $\GL_{n+1}(\C)$ on $\Pi_{\dd,n}$ preserves $\Sigma_{\dd,n}$.

Take $(f^0_1,\ldots,f^0_k)\in \Pi_{\dd,n}\setminus\Sigma_{\dd,n}$; denote by $\bb_i^{\dd,n}$ the pullback of $\aaa_i^{\dd,n}$ under
the corresponding orbit map $\GL_{n+1}(\C)\to\Pi_{\dd,n}\setminus\Sigma_{\dd,n}$. We obviously have
\label{b}$$\bb_i^{\dd,n}=\m_i^{\dd,n}\cc_i^{n+1}+o(\cc_i^{n+1}),$$
where \label{m}$\m_i^{\dd,n}$ are integers (notice that these integers do not depend on the choice of $(f^0_1,\ldots,f^0_k)\in \Pi_{\dd,n}\setminus\Sigma_{\dd,n}$).
One of our main tasks in the sequel will be obtaining
%estimates (and in some cases,
explicit expressions
%)
for $\m_i^{\dd,n}$.

\subsubsection{Miscellany}
For any $i=1,\ldots,k$ define the suspension map \label{SSS}$S^{\dd,n}_i:\Pi_{\dd,n}\to\Pi_{\dd,n+1}$ by the formula
\begin{equation}\label{suspension}
(f_1,\ldots,f_k)\mapsto(f_1,\ldots,f_i+x_{n+1}^{d_i},\ldots,f_k)
\end{equation}
(here $f_1,\ldots,f_k$ are polynomials in $x_0,\ldots,x_n$).

Recall also the following identity (sometimes called the Euler formula):
\begin{equation}\label{euler}f(x_0,\ldots ,x_n)=\frac{1}{d}\sum_{j=0}^n x_j\frac{\partial f}{\partial x_j}(x_0,\ldots ,x_n),\end{equation}
where $f$ is a homogeneous complex polynomial of degree $d$.

\section{The action of $\GL_{n+1}(\C)$ on $\Pi_{\dd,n}\setminus\Sigma_{\dd,n}$}\label{sec3}
This section and the following one are devoted to the calculation of $\m_i^{\dd,n}$ (and hence, the details get somewhat technical at times).
\subsection{The calculation of $\m_{n+1}^{\dd,n}$}
Take $(f^0_1,\ldots,f^0_k)\not\in\Sigma_{\dd,n}$.
Let $F_{\dd,n}$ be the mapping $\C^{n+1}\to\Mat_{n+1,k}(\C)\times\C^{n+1}$ defined by
\begin{equation}\label{defndn}x\mapsto(df^0_1|_x,\ldots,df^0_k|_x,x).\end{equation}
Due to (\ref{euler}), $F_{\dd,n}(x)\in\W_{k,n}$, iff $x=0$. Define $N(\dd,n)$ by the formula
\label{defdefndn}$$F_{\dd,n}^*(\lk_{[\W_{k,n}],\Mat_{n+1,k}(\C)\times\C^{n+1}})=N(\dd,n)\ee_{n+1}.$$

%The following proposition will allow us to deduce the division theorem \ref{main}, without actually computing $N(\dd, k)$.
In the sequel we give (proposition \ref{krainie} and lemma \ref{ladders}) explicit formulae for $N(\dd,n)$; however we present here one basic property of these numbers, since, on the one hand, it follows directly from the definition, and on the other hand, it is all we shall need for the proof of theorem \ref{main}. 

\begin{Prop}\label{Ndn}
We have $N(\dd, n)>1$, unless $\dd=(2)$.
\end{Prop}

{\bf Proof.} If $k>1$, the degree of $\W_{k,n}$ is $>1$, and zero belongs to the singular locus of $\W_{k,n}$. Hence, the intersection multiplicity of $\W_{k,n}$ and the image of $F_{\dd,n}$ at 0 is $>1$.

If $k=1$, the sequence $\dd$ contains just one element, $\dd=(d)$. The degree of $\W_{1,n}$ is 1, and
$N((d),n)$ is equal to the degree of the mapping $\C^{n+1}\ni x\mapsto df^0|_x,f^0\in\Pi_{(d),n}\setminus\Sigma_{(d),n}$, which is $(d-1)^{n+1}>1$,
%independent of $f^0_1\in\Pi_{(d),n}\setminus\Sigma_{(d),n}$. By taking $f^0_1=\sum_{i=0}^n x_i^d$, we see that the degree of $F_{(d),n}$ is $(d-1)^{n+1}>1$,
unless $\dd=(2)$.$\clubsuit$

\begin{lemma}\label{highnonempty}
We have $\bb_{n+1}^{\dd,n}=(N(\dd,n)+(-1)^{n+k+1})\cc_{n+1}^{n+1}$, and hence, $\m_{n+1}^{\dd,n}=N(\dd,n)+(-1)^{n+k+1}$.
\end{lemma}

{\bf Proof.}
%We proceed as in the proof of lemma \ref{emptyhigh}.
Take $(f^0_1,\ldots,f^0_k)\not\in\Sigma_{\dd,n}$, and set $x_0=(0,\ldots,0,1)^T$.

Recall that the class $\aaa_{n+1}^{\dd,n}$ is the linking number with
the variety $$\{(f_1,\ldots,f_k)|f_1(x_0)=\cdots=f_k(x_0)=0,df_1|_{x_0},\ldots,df_k|_{x_0}\mbox{are linearly dependent}\}.$$
Due to (\ref{euler}), this variety is the preimage of $\W_{k,n}$ under
$$(f_1,\ldots,f_k)\mapsto (df_1|_{x_0},\ldots,df_k|_{x_0},x_0).$$
This mapping $\Pi_{\dd,n}\to\Mat_{n+1,k}(\C)\times\C^{n+1}$ can be included into the following commutative diagram
$$\xymatrix{\GL_{n+1}(\C)\ar[dr]\ar[d]& \\\Pi_{\dd,n}\ar[r]&\Mat_{n+1,k}(\C)\times\C^{n+1}}$$
where the vertical arrow is the orbit map, and the diagonal one is the map
\begin{equation}\label{map}A\mapsto (A^Tdf^0_1|_{Ax_0},\ldots,A^Tdf^0_k|_{Ax_0},x_0).\end{equation}

Denote by $\alpha$ the linking number
$$\lk_{[\W_{k,n}\cap(\Mat_{n+1,k}(\C)\times\{x_0\})],\Mat_{n+1,k}(\C)\times\{x_0\}}.$$
Applying corollary \ref{linear}, we see that $\bb_{n+1}^{\dd,n}$ is the pullback of $\alpha$ under (\ref{map}).

Set ${\cal E}=\Mat_{n+1,k}(\C)\times(\C^{n+1}\setminus\{0\})$, and set $V={\cal E}\cap\W_{k,n}$. Denote by $\beta$ the restriction of $\lk_{[\W_{k,n}],\Mat_{n+1,k}(\C)\times\C^{n+1}}$
to ${\cal E}\setminus V$. Since $V$ if fibered over $\C^{n+1}\setminus\{0\}$, the intersection of
$V$ and $\Mat_{n+1,k}(\C)\times\{x\}$ is transversal for any $x\in\C^{n+1}\setminus\{0\}$, which implies (due to the second assertion of proposition \ref{lk}) that the restriction of $\beta$ to
$$(\Mat_{n+1,k}(\C)\times\{x_0\})\setminus\W_{k,n}=(\Mat_{n+1,k}(\C)\times\{x_0\})\setminus V$$ is $\alpha$.

The group $\GL_{n+1}(\C)$ acts on ${\cal E}$ by the formula
$$A\cdot(B,x)=((A^T)^{-1}B,Ax).$$ This action preserves $V$; indeed, if $x^T B=0$, then $(Ax)^T(A^T)^{-1}B=0$. Let us compute the pullback of $\beta$ under the
action map $$\GL_{n+1}(\C)\times ({\cal E}\setminus V)\to {\cal E}\setminus V.$$ The nontrivial cohomology groups of ${\cal E}\setminus V$ start in dimension $2n+1$, which implies that
the pullback of $\beta$ has the form $\gamma\otimes 1+1\otimes\beta$, where $\gamma\in H^{2n+1}(\GL_{n+1}(\C))$.

\begin{lemma}\label{gamma}
We have $\gamma=(-1)^{n+k}\cc^{n+1}_{n+1}$.
\end{lemma}

We shall prove this lemma a little later.

The mapping (\ref{map}) can be represented as the composition
\begin{equation}
\label{comp}\GL_{n+1}(\C)\to\xymatrix{\GL_{n+1}(\C)\times(\C^{n+1}\setminus\{0\})\ar[r]^-{Id\times F_{\dd,n}}&
\GL_{n+1}(\C)\times ({\cal E}\setminus V)}\to {\cal E}\setminus V\end{equation}
(where the first map is $A\mapsto (A^{-1},\mbox{the last column of $A$})$, and
the third one is the action map). Applying lemma \ref{gamma} and the second assertion of proposition \ref{glm}, we obtain that
the pullback of $\beta$ under \ref{comp} is $-\gamma+N(\dd,n)\cc^{n+1}_{n+1}=(N(\dd,n)+(-1)^{n+k+1})\cc_{n+1}^{n+1}$, which implies lemma
\ref{highnonempty}.$\clubsuit$

\subsubsection{Proof of lemma \ref{gamma}}
We keep the notation of the proof of lemma \ref{highnonempty}.

%Let us calculate the pullback of $\beta$ under the analogue of (\ref{map}) for linear polynomials.
Define the (linear) polynomial $g_i,i=1,\ldots,k$, by $g_i=x_{n-i+1}$, and set
$$F_\mathrm{lin}(x)=(dg_1|_x,\ldots,dg_k|_x,x),$$
$$s(A)=(A^Tdg^0_1|_{Ax_0},\ldots,A^Tdg^0_k|_{Ax_0},x_0)=(\mbox{the last $k$ columns of $A^T$},x_0)\eqno (\ref{map}')$$
(recall that $x_0=(0,\ldots,0,1)^T$). The mapping $s$ can be represented as the composition
$$\GL_{n+1}(\C)\to\GL_{n+1}(\C)\times(\C^{n+1}\setminus\{0\})\to
\GL_{n+1}(\C)\times ({\cal E}\setminus V)\to {\cal E}\setminus V,\eqno (\ref{comp}')$$
where the first and the last arrows are the same as in (\ref{comp}), and the middle one is $Id\times F_{\mathrm{lin}}$.
The image of $F_\mathrm{lin}$ does not meet $\W_{k,n}$, hence, pullback of $\beta$ under (\ref{map}${}'$) is $-\gamma$. We shall now calculate
this pullback directly, which will complete the proof of lemma \ref{gamma}.

Recall that above we have defined
$$W_{k,n}=\{A\in\Mat_{n+1,k}(\C)\mid \rk A<k\},$$
$$X_{k,n}=\{A\in\Mat_{n+1,k}(\C)\mid \rk A<k,\mbox{the last row of $A$ is zero}\}.$$
Notice that $X_{k,n}\times\{x_0\}=\W_{k,n}\cap(\Mat_{n+1,k}(\C)\times\{x_0\})$.
Let $\xi_0,\xi_1$ and $\xi_2$ be the vector bundles on $\C P^{k-1}$ such that
$$\tot (\xi_0)=\{((y_1,\ldots,y_k),(z_1:\ldots:z_k))\in\C^{k}\times\C P^{k-1}\mid \sum y_iz_i=0\},$$
$$\tot (\xi_1)=\{(A,z)\in\Mat_{n+1,k}(\C)\times\C P^{k-1}\mid Az=0\},$$
$$\tot (\xi_2)=\{(A,z)\in\Mat_{n+1,k}(\C)\times\C P^{k-1}\mid Az=0,\mbox{the last row of $A$ is zero}\}.$$

Denote by $a$ the image of $\ee_{n+1}\in H^{2n+1}(\C^{n+1}\setminus\{0\})$ under the map
$$(\Mat_{n+1,k}(\C)\setminus W_{k,n})\ni A\mapsto\mbox{the last column of $A$},$$ and let $b$ be
the maximal power of the the canonical generator of $H^*(\C P^{k-1})$.
Set $b'$ to be the image of $b$ under the isomorphism $H^*(\C P^{k-1})\to H^*(\tot(\xi_2))$.

\begin{Prop}\label{chernclass}
We have $c_{k-1}(\xi_0)=e(\xi_0)=(-1)^{k+1}b$.
\end{Prop}

{\bf Proof.} The direct sum of $\xi_0$ and the cotautological bundle is isomorphic (as a topological vector bundle) to the trivial rank $k$ bundle on $\C P^{k-1}$.$\clubsuit$

\begin{Prop}\label{restr}
The restriction of $\lk_{[X_{k,n}],\Mat_{n+1,k}(\C)}$ to $\Mat_{n+1,k}(\C)\setminus W_{k,n}$ is $(-1)^{k+1}a$.\end{Prop}

{\bf Proof of proposition \ref{restr}.} Let us first calculate the image of $[X_{k,n}]$ in $\bar H_*(W_{k,n})$.
Set $\xi$ to be the pullback of $\xi_0$ to $\tot\xi_2$. Notice that $\tot\xi$ can be naturally identified with $\tot\xi_1$.
%The space $\tot(\xi_1)$ may be viewed as the total space of the pullback $\xi$
%of $\xi_0$
%to $\tot(\xi_2)$.

Let $F$ be the fibre of $\xi_2$ over $(0:\cdots :0:1)$; clearly, $F$ is Poincar\'e dual to $b'$ in $\tot(\xi_2)$, hence, by propositions \ref{zerosec} and \ref{chernclass},
the image of $[\tot (\xi_2)]$
in $\bar H_*(\tot(\xi_1))$ is $(-1)^{k+1}[\tot(\xi|_F)]=(-1)^{k+1}[\mbox{the fibre of $\xi_1$ over $(0:\cdots:0:1)$}]$.
The diagram
\begin{equation*}
\xymatrix{\tot(\xi_2)\ar[r]\ar[d]&X_{k,n}\ar[d]\\ \tot(\xi_1)\ar[r]&W_{k,n}}\end{equation*} implies then that the image of $[X_{k,n}]$ in $\bar H_*(W_{k,n})$ is $(-1)^{k+1}(\mbox{the image of $[U]$})$, where $U$ is
is the union of all matrices $\in\Mat_{n+1,k}(\C)$ with zero last column.

%By the third assertion of proposition \ref{lk}, the
The restrictions of $\lk_{[X_{k,n}],\Mat_{n+1,k}(\C)}$ and $(-1)^{k+1}\lk_{[U],\Mat_{n+1,k}(\C)}$
to $\Mat_{n+1,k}(\C)\setminus W_{k,n}$ coincide (cf. diagram (\ref{diagxy})). Applying corollary \ref{linear}, we obtain that the restriction of
$\lk_{[U],\Mat_{n+1,k}(\C)}$ to $\Mat_{n+1,k}(\C)\setminus W_{k,n}$ is $a$, and
proposition \ref{restr} follows.$\clubsuit$

Proposition \ref{restr} and the first assertion of proposition \ref{glm} imply that the pullback of $\beta$ under
(\ref{map}${}'$) is $(-1)^{n+k+1}\cc^{n+1}_{n+1}=-\gamma$. Lemma \ref{gamma} is now proven.$\clubsuit$

\subsection{The calculation of the remaining $\m_i^{\dd,n}$.}
%Recall that we suppose $k<n+1$.
\label{dd'}If $k>1$, set $\dd'=(d_2,\ldots,d_k)$.
Suppose we have proven that $\bb_i^{\dd,n-1}=\m_i^{\dd,n-1}\cc_i^n$ for $k<n+1,i\leq n$ and $\bb_i^{\dd',n-1}=m_i^{\dd,n-1}\cc_i^n.$
We shall express here
$\m_i^{\dd,n}$ in terms of $\m_i^{\dd,n-1}$ and
$\m_i^{\dd',n-1}$.

Consider first the ``generic'' case $k<n+1$. Take an element $(f_1,\ldots,f_k)\in\Pi_{\dd,n-1}$, and take a point
$(0:\cdots :0:x_{i-1}:\cdots :x_{n-1})\in\C P^{n-1}$. There exists $x_n$ 
such that
$\Sing(f_1+x_n^{d_1},\ldots,f_k)\ni (0:\cdots :0:x_{i-1}:\cdots :x_{n-1}:x_n)$,
iff
$$(0:\cdots :0:x_{i-1}:\cdots :x_{n-1})\in\Sing (f_1,\ldots,f_k)\cup\Sing (f_2,\ldots,f_k).$$ This implies that the preimage
of $V_{\dd,n,\C P^{n+1-i}}$ under the suspension map $S_1^{\dd,n-1}$
is the union of $V_{\dd,n-1,\C P^{n-i}}$ and
$$\{(f_1,\ldots,f_k)\mid (f_2,\ldots,f_k)\in V_{\dd',n-1,\C P^{n-i}}\}$$ (if $k=1$, the second set is empty).
Denote by \label{deltagamma}$\delta_i^{\dd,n}$ and $\gamma_i^{\dd,n}$ the corresponding
multiplicities (set $\gamma_i^{\dd,n}=0$, if $k=1$).
%and if $k=n+1$, we set $\delta_i^{\dd,n}=0$).
Identify $\GL_n(\C)$ with the image of the canonical inclusion $\GL_n(\C)\to\GL_{n+1}(\C)$; under this identification the map $S_1^{\dd,n}$ becomes
$\GL_n(\C)$-equivariant, and we obtain (using the second assertion of proposition \ref{lk}) that the restriction of $\bb_{i}^{\dd,n},i=1,\ldots,n,$ to $\GL_n(\C)$
is $(\gamma_i^{\dd,n}\m_i^{\dd',n-1}
+\delta_i^{\dd,n}\m_i^{\dd,n-1})\cc_i^n$.

Suppose now $k=n+1>1$ (the case $k=1,n=0$ corresponding to
empty hypersurfaces in $\C P^0$ has in fact already been considered in the previous subsection).
Set $f_1=x_n^{d_1}$.
%Suppose that $n>0$ and that we have calculated $\m_1^{\dd',n-1},\ldots,\m_{n-1}^{\dd',n-1}$.
We replace the suspension map $S_1^{\dd,n-1}$
%(which, in fact, is not defined in the case $k=n+1$, since the lenght of $\dd$ is $k=n+1>(n-1)+1$)
by the map \begin{equation}\label{suspempty}(f_2,\ldots,f_k)\mapsto (f_1,f_2,\ldots,f_k).\end{equation}
A point $(0:\cdots 0:x_{i-1}:\cdots:x_{n-1})\in\Sing(f_2,\ldots,f_k)$ (i.e.,
the polynomials
$f_2,\ldots,f_k$ have a common zero at $(0,\ldots,0,x_{i-1},\ldots,x_{n-1})\in\C^n$),
iff there exists $x_n$ such that $f_1,f_2,\ldots,f_k$ have a common zero at $(0,\ldots,0,x_{i-1},\ldots,x_{n-1},x_n)\in\C^{n+1}$
(in fact, such $x_n$ is necessarily zero).
This implies that the intersection of the image of (\ref{suspempty})
%\begin{equation}\label{suspempty}(f_2,\ldots,f_k)\mapsto (f_1,f_2,\ldots,f_k)\end{equation}
and $V_{\dd,n,\C P^{n+1-i}}$
%under the map \begin{equation}\label{suspempty}(f_2,\ldots,f_k)\mapsto (f_1,f_2,\ldots,f_k)\end{equation} is $V_{\dd',n-1,\C P^{n-i}}$.
is $$\{(f_1,\ldots,f_k)\mid (f_2,\ldots,f_k)\in V_{\dd',n-1,\C P^{n-i}}\}$$
Denote by $\gamma_i^{\dd,n}$ the corresponding multiplicity.
If we identify as above $\GL_n(\C)$ with the image of the canonical inclusion $\GL_n(\C)\to\GL_{n+1}(\C)$, the map (\ref{suspempty}) becomes
$\GL_n(\C)$-equivariant, hence  $\bb_{i}^{\dd,n},i=1,\ldots,n,$ restricted to $\GL_n(\C)$
is $\gamma_i^{\dd,n}\m_i^{\dd',n-1}\cc_i^n$.
%The proof of the following lemma is then completed analogously
%to the proof of lemma \ref{xxx}.

We summarise the results of this subsection in the following lemma (where $\delta_i^{\dd,n}$ is set to be zero in the case $k=n+1$).

\begin{lemma}\label{prevnonempty}
Suppose that $n>0$. For any $i=1,\ldots,n,$ there exist integers $\gamma_i^{\dd,n},\delta_i^{\dd,n}\geq 0$ at least one on which is positive and such that
$\bb_i^{\dd,n}=(\gamma_i^{\dd,n}\m_i^{\dd',n-1}
+\delta_i^{\dd,n}\m_i^{\dd,n-1})\cc_i^{n+1},$ and hence,
$\m_i^{\dd,n}=\gamma_i^{\dd,n}\m_i^{\dd',n-1}
+\delta_i^{\dd,n}\m_i^{\dd,n-1}$.
\end{lemma}
$\clubsuit$

\section{Fixing the coefficients}\label{sec4}
In the previous section we have introduced a lot of various coefficients.
Here we convert them into numbers.
\subsection{The calculation of $N(\dd,n)$}
Let us first consider the ``extreme'' cases $k=1$ and $k=n+1$.
\begin{Prop}\label{krainie}
\begin{enumerate}
\item For any $d\geq 2$ we have $N((d),n)=(d-1)^{n+1}$.
\item If $k=n+1$, then $N(\dd,n)=\prod_{i=1}^k d_i-1$.
\end{enumerate}
\end{Prop}
(Notice that in the most extreme case $k=n+1=1$ both formulae coincide.)

{\bf Proof.}
%In the case $k=1$ the variety $\W_{k,n}$ is just a vector subspace of $\Mat_{n+1,k}(\C)\times\C^{n+1}$, and hence, $N((d),n)$ is the degree
%of the maping $\C^{n+1}\ni x\mapsto df_0|_x$, where $f_0$ is any polynomial defining a smooth hypersurface. By setting $f_0=\sum x_i^d$ we obtain
%We have already seen in the proof of proposition \ref{Ndn} that 
The first one of these formulae has already been obtained in the proof of proposition \ref{Ndn}, so it remains to prove the second one. Suppose $k=n+1$.
Notice that in this case we can compute $\m_{n+1}^{\dd,n}$ directly, i.e., without using lemma~\ref{highnonempty}. Indeed, set $x_0=(0:\cdots:0:1)$; clearly, $x_0\in\Sing(f_1,\ldots,f_{n+1})$, iff $f_1(x_0)=\cdots=f_{n+1}(x_0)=0$. Hence, $V_{\dd,n,\{x_0\}}$ is the preimage of $0\in\C^{n+1}$ under the mapping $$(f_1,\ldots,f_{n+1})\mapsto (f_1(x_0),\ldots,f_{n+1}(x_0))^T.$$ Hence, the class $\bb_{n+1}^{\dd,n}$ is equal to $\cc_{n+1}^{n+1}$ times the degree $\lambda$ of the mapping $$x\mapsto (f^0_1(x),\ldots,f^0_{n+1}(x))^T,$$ where $(f^0_1,\ldots,f^0_{n+1})\in\Pi_{\dd,n}\setminus\Sigma_{\dd,n}$ (cf. the proof of lemma \ref{highnonempty}). By setting $f^0_i=x_{i-1}^{d_i}$, we see that $\lambda=\prod_{i=1}^{n+1}d_i$. The required formula for $N(\dd,n)$ follows now from lemma \ref{highnonempty}.$\clubsuit$
%The ideal $I(\W_{n+1,n})$ is generated by the polynomial $(A,x)\mapsto\det A$ (where $A\in\Mat_{n+1,n+1}(\C), x\in\C^{n+1}$) and $n+1$ quadratic %polynomials obtained from the relation $A^Tx=0$.
%By setting in (\ref{defndn}) $f_i^0=x_{i-1}^{d_i},i=1,\ldots,n+1,$ we see that the pullback of $I(\W_{n+1,n})$ under $F_{\dd,n}$ is
%generated by $x_{i-1}^{d_i},i=1,\ldots ,n+1,$ and by $\prod_{i=0}^nx_{i-1}^{d_i-1}$. The proof is now completed using the third assertion of proposition \ref{lk}.$\clubsuit$
%
%Before considering the ``generic'' case $1<k<n+1$, let us make the following obvious observation, which we shall use a few times in this section.
%
%\begin{Prop}
%Let $F_1:\C^m\to\C^N,F_2:\C\to\C^2$ be polynomial mappings.
%Suppose that $V\subset\C^N$ and $W\subset\C^N\times\C^2=\C^{N+2}$ are irreducible affine varieties
%such that $\codim_{\C^N}V=\codim_{\C^{N+2}}W=m$ and $W\subset V\times\C^2$. Denote by $W_x,x\in V$ (the projection to $\C^2$ of) the intersection
%$W\cap(\{x\}\times\C^2)$. Set $F=F_1\times F_2$.
%Suppose that for any $x\in V$ we have $\codim_\C F_2^{-1}(I(W_x))=l>0$. Then $\codim_\C F^{-1}(I(W))=l\codim_\C F_1^{-1}(I(V))$.
%\end{Prop}
%$\clubsuit$

\begin{Prop}\label{recndn}
Suppose that $1<k<n+1$.
%, and set $\dd'=(d_2,\ldots,d_k)$.
Then we have
$N(\dd,n)=(d_1-1)N(\dd,n-1)+d_1N(\dd',n-1)$.
\end{Prop}
%(Recall that $\dd'$ is by definition equal to $(d_2,\ldots,d_k)$.)

{\bf Proof.} Take $(g_1,\ldots,g_k)\in\Pi_{\dd,n-1}\setminus\Sigma_{\dd,n}$ such that $(g_2,\ldots,g_k)\in\Pi_{\dd',n-1}\setminus\Sigma_{\dd',n-1}$, and set in
(\ref{defndn})
$f^0_1=g_1+x_n^{d_1},f^0_i=g_i,i=2,\ldots,k$. The image of $F_{\dd,n}$ will be then contained in the vector subspace of $\mathbf{E}\subset\Mat_{n+1,k}\times\C^{n+1}$
defined by (\ref{matrix}). Applying proposition \ref{matrixlemma} and the second assertion of proposition \ref{lk}, we have
$$F_{\dd,k}^*(\lk_{[\W_{\dd,n}],\Mat_{n+1,k}\times\C^{n+1}})=F_{\dd,k}^*(\lk_{[\X_{k,n}],\mathbf{E}}+\lk_{[\Y_{k,n}],\mathbf{E}}).$$

Let us calculate $F_{\dd,k}^*(\lk_{[\X_{k,n}],\mathbf{E}})\in H^{2n+1}(\C^{n+1}\setminus\{0\}$.
Consider the isomorphism $R:\mathbf{E}\to (\Mat_{n,k}\times\C^n)\times\C^2$ defined as follows:
take an element $(A,x)\in \mathbf{E}$, where $$A=(a_{i,j})_{\genfrac{}{}{0pt}{}{0\leq i\leq n}{1\leq j\leq k}},x=(x_0,\ldots,x_n)^T,$$ to
$((A',x'),a_{n,1},x_n)$, where $$A'=(a_{i,j})_{\genfrac{}{}{0pt}{}{0\leq i\leq n-1}{1\leq j\leq k}},x'=(x_0,\ldots,x_{n-1})^T$$
\begin{figure}\centering
\epsfbox{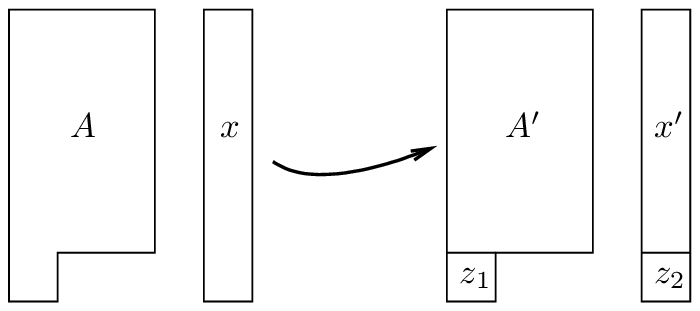}
\caption{The isomorphism $R$}\label{isor}
\end{figure}
(see figure \ref{isor}).

Under this identification, the variety $\X_{k,n}$ is taken to $\W_{k,n-1}\times\{(z_1,z_2)\mid z_1=0\},$ and the mapping $F_{\dd,n}$
is written as
$$F_{\dd,n-1}\times (\mbox{the map $(x_0,\ldots,x_n)^T\mapsto (d_1 x_n^{d_1-1},x_n)$}),$$
%$$(x_0,\ldots,x_n)\mapsto (F_{\dd,n-1}(x_0,\ldots,x_{n-1}),(d_1 x_n^{d_1-1},x_n)),$$
where $F_{\dd,n-1}$ is obtained by setting in $f^0_i=g_i,i=1,\ldots,k,$ in (\ref{defndn}). This proves that
$$F_{\dd,k}^*(\lk_{[\X_{k,n}],\mathbf{E}})=(d_1-1)N(\dd,n-1)\ee_{n+1}.$$

In order to complete the proof of the proposition, it is sufficient to show that
\begin{equation}\label{xxy}F_{\dd,k}^*(\lk_{[\Y_{k,n}],\mathbf{E}})=d_1N(\dd',n-1)\ee_{n+1}.\end{equation} Indeed, consider the isomorphism
$T:\mathbf{E}\to\C^n\times (\Mat_{n,k-1}(\C)\times\C^n)\times\C^2$ defined as follows: we take
$(A,x)\in \mathbf{E}$, where $$A=(a_{i,j})_{\genfrac{}{}{0pt}{}{0\leq i\leq n}{1\leq j\leq k}},x=(x_0,\ldots,x_n)^T,$$ to
$(v,(A',x'),a_{n,1},x_n)$, where $$A'=(a_{i,j})_{\genfrac{}{}{0pt}{}{0\leq i\leq n-1}{2\leq j\leq k}},x'=(x_0,\ldots,x_{n-1})^T,v=(a_{1,1},\ldots,a_{n,1})^T$$
\begin{figure}\centering
\epsfbox{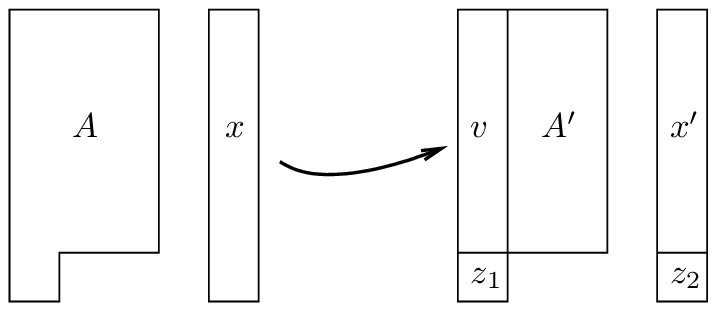}
\caption{The isomorphism $T$}\label{isot}
\end{figure}
(see figure \ref{isot}).

Under this identification, $\Y_{\dd,n}$ is transformed into
\begin{equation}\label{TYkn}\{(v,(A',x'),(z_1,z_2))\mid (A',x')\in\W_{k-1,n-1},v^T\cdot x'+z_1z_2=0\},\end{equation}
and the mapping
$F_{\dd,n}$ is becomes $$\left(\mbox{the map $\C^{n+1}\ni\begin{pmatrix}x\\ 0\end{pmatrix}\mapsto dg_1|_x\in\C^n$}\right)\times F_{\dd',n-1}\times(\mbox{the map $(x_0,\ldots,x_n)^T\mapsto (d_1 x_n^{d_1-1},x_n)$}),$$
where $F_{\dd',n-1}$ is obtained by setting $f^0_i=g_{i+1},i=1,\ldots,k-1$ in (\ref{defndn}). The variety
(\ref{TYkn}) projects onto $\C^n\times\W_{k-1,n-1}$, the preimage of a point being a curve $\subset\C^2$ of the form $\{(z_1,z_2)\mid z_1z_2=a\}$.
Hence, the intersection multiplicity of the image of $F_{\dd,n}$ and $\Y_{k,n}$ is equal to $d_1$ times the intersection multiplicity
of the image of $F_{\dd',n-1}$ and $\W_{k-1,n-1}$, which proves (\ref{xxy}).$\clubsuit$
%We can identify $\mathbf{E}$ with $(\Mat_{n,k}\times\C^n)\times\C^2$
%by taking $(A,x),A=(a_{i,j}),x=(x_0,\ldots,x_n)$, to $((A',x'),a_{n+1,1},x_n)$, where $A'$ (respectively, $x'$) is obtained from $A$
%(respectively, $x$) by forgetting the last row (respectively, the last coordinate).
%so that $X_{k,n}$ will become
%$$\W_{k,n-1}\times\{(z_1,z_2)\mid z_1=0\},$$ and the mapping $F_{\dd,n}$ will be written
%$$(x_0,\ldots,x_n)=(F_{\dd,n-1}(x_0,\ldots,x_{n-1}),(d_1 x_n^{d_1-1},x_n)).$$

Let us present the results of this subsection in a (more or less) compact way.
\begin{lemma}\label{ladders}
We have for $k>1$
\begin{multline*}
N(\dd,n)=(d_1-1)^{n-k+2}\\
+\sum_{i=2}^k\left.\frac{d_{1}\cdots d_{i-1}(d_i-1)}{(n-k+1)!}\frac{d^{n-k+1}}{dt^{n-k+1}}\right|_{t=0}\frac{1}{(1-(d_1-1)t)\cdots(1-(d_i-1)t)}.
\end{multline*}
\end{lemma}
This expression does not look symmetric in $d_1,\ldots,d_k$ but in fact it is.

{\bf Proof.} We can arrange the $N(\dd,n)$'s into a table the following one (represented here for $n=5,k=3$):
$$\begin{array}{c|ccc}
%\cline{2-4}
2&(d_3-1)^4&N((d_2,d_3),4)&N((d_1,d_2,d_3),5)\\
%\cline{2-4}
1&(d_3-1)^3&N((d_2,d_3),3)&N((d_1,d_2,d_3),4)\\
%\cline{2-4}
0&(d_3-1)^2&N((d_2,d_3),2)&N((d_1,d_2,d_3),3)\\
%\cline{2-4}
-1&d_3-1&d_2d_3-1&d_1d_2d_3-1\\
%\cline{2-4}
-2&1&1&1\\
\hline
&1 & 2 & 3
\end{array}$$
(Here the $x$-coordinate is $k$, the $y$-coordinate is $n-k$, and the bottom line of 1's is added for formal reasons.)

A path connecting two boxes of this table will be called a {\it staircase}, if it goes only downwards or to the left.
A {\it segment} is a staircase that joins two neighbouring boxes.
Let us associate the {\it weight} $d_{k-a+1}$, respectively, $d_{k-a+1}-1$, to the segment joining the boxes $(a,b)$ and $(a-1,b)$,
respectively, the boxes $(a,b)$ and $(a,b-1)$; the weight of a staircase is set to be the product of the weights of all of its segments.
Due to propositions \ref{krainie} and \ref{recndn}, $N(\dd,n)$ is equal to the sum of the weights of all staircases
that descend from the box containing $N(\dd,n)$ to the bottom line of 1's and whose last segment is vertical.$\clubsuit$

\subsection{The multiplicities $\gamma_i^{\dd,n}$ and $\delta_i^{\dd,n}$}
\begin{lemma}\label{delta} If $k<n+1$, then $\delta_i^{\dd,n}=d_1-1$.
\end{lemma}

{\bf Proof.} We can obtain $\delta_i^{\dd,n}$ as follows. Let $U\subset\Pi_{\dd,n-1}$ be a small $i$-dimensional ball that
intersects $V_{\dd,n-1,\C P^{n-i}}$ transversally at a smooth point $\{(f_1^0,\ldots,f_k^0)\}$
(the subspace $\C P^{n-i}$ in the definition of $V_{\dd,n-1,\C P^{n-i}}$
is taken to be
\begin{equation}\label{cpn-i}
\{(0:\cdots :0:x_{i-1}:\cdots :x_{n-1})\in\C P^{n-1}\},\end{equation} as in the proof of lemma \ref{prevnonempty}).

Assume that $\Sing(f_1^0,\ldots,f_k^0)\cap\C P^{n-i}=\{(0:\cdots :0:1)\}$.
Moreover, $V_{\dd,n-1,\C P^{n-i}}$ and $\Pi_{(d_1),n-1,\C P^{n-i}}\times V_{\dd',n-1,\C P^{n-i}}$ are distinct irreducible subvarieties of $\Pi_{\dd,n-1}$ of the same dimension, hence, we can also assume $U\cap(\Pi_{(d_1),n-1,\C P^{n-i}}\times V_{\dd',n-1,\C P^{n-i}})=\varnothing$.
%that $(f_2,\ldots,f_k)\not\in V_{\dd',n-1,\C P^{n-i}}$
%for any $(f_2,\ldots,f_k)$ such that $(f_1,f_2,\ldots,f_k)\in U$.
Clearly, $\delta_i^{\dd,n}$ is equal to the intersection multiplicity $\mu$ of $U'=S_1^{\dd,n-1}(U)$ and $V_{\dd,n,\C P^{n-i+1}}$ at
$S_1^{\dd,n-1}(f_1^0,\ldots,f_k^0)=(f_1+x_n^{d_1},f_2,\ldots,f_k)$ (once again, we take the subspace
\begin{equation}\label{cpn-i+1}
\C P^{n-i+1}=\{(0:\cdots :0:x_{i-1}:\cdots :x_{n})\in\C P^{n}\}\end{equation} used in the proof of lemma \ref{prevnonempty} as the subspace
$\C P^{n-i+1}$ in the definition of $V_{\dd,n,\C P^{n-i+1}}$). Set
%Set $U_1$ (respectively, $U_2$) to be the following neighbourhood of $(0:\cdots:0:1)$ in $\C P^{n-i}$
%(respectively, of $(0:\cdots:0:1:0)$ in $\C P^{n-i+1}$):
%$$U_1=\{(0:\cdots :0:z_1:\cdots:z_{n-i}:1:0)\},$$
$$U_1=\{(0:\cdots :0:z_1:\cdots:z_{n-i}:1:z_{n-i+1})\}.$$
For any $x=(0:\cdots :0:z_1:\cdots:z_{n-i}:1:z_{n-i+1})\in U_1$ set $\tilde x$ to be the lifting $(0,\ldots ,0,z_1,\ldots,z_{n-i},1,z_{n-i+1})^T$ of $x$.
Finally, set
%$\tilde V_{\dd,n-1,\C P^{n-i}}$ and
$\tilde V_{\dd,n,\C P^{n-i+1}}$ to be the natural smooth
resolution
%s of $V_{\dd,n-1,\C P^{n-i}}$ and 
$V_{\dd,n,\C P^{n-i+1}}$,
%respectively,
i.e.,
%$$\tilde V_{\dd,n-1,\C P^{n-i}}=\{((f_1,\ldots,f_k),x)\in\Pi_{\dd,n-1}\times\C P^{n-i}\mid x\in\Sing(f_1\ldots,f_k)\},$$
$$\tilde V_{\dd,n,\C P^{n-i+1}}=\{((f_1,\ldots,f_k),x)\in\Pi_{\dd,n}\times\C P^{n-i+1}\mid x\in\Sing(f_1\ldots,f_k)\}.$$
%and denote by $p_1$ and $p_2$ the obvious projections $\tilde V_{\dd,n-1,\C P^{n-i}}\to V_{\dd,n-1,\C P^{n-i}}$ and
%$\tilde V_{\dd,n,\C P^{n-i+1}}\to V_{\dd,n,\C P^{n-i+1}}$.

Obviously, $\mu$ is equal to the intersection multiplicity of $U'\times U_1$ and $\tilde V_{\dd,n,\C P^{n-i+1}}$ at
$$((f_1+x_n^{d_1},f_2,\ldots,f_k),(0:\cdots:0:1:0)),$$ which is equal to the the intersection multiplicity of $\W_{\dd,n}$ and the image
of the mapping $F:U\times U_1\to\Mat_{n+1,k}(\C)\times\C^{n+1}$ given by
%that takes $$((f_1,\ldots,f_k),x),x=(0:\cdots :0:z_1:\cdots:z_{n-i}:1:z_{n-i+1}),$$ to
$$((f_1,\ldots,f_k),x)\mapsto (df_1|_{\tilde x},\ldots,df_k|_{\tilde x},\tilde x)$$
%
%(indeed, the restriction of $\tilde V_{\dd,n,\C P^{n-i+1}}$ to $U_1$ is the preimage of the bundle
%$$\W_{k,n}\cap (\Mat_{n+1,k}(\C)\times (\C^{n+1}\setminus\{0\}))\to\C^{n+1}\setminus\{0\}$$ under the mapping
%$((f_1,\ldots,f_k),x)\mapsto (df_1|_{\tilde x},\ldots,df_k|_{\tilde x},\tilde x)$).

It can be readily seen that $F$ is a local embedding at $((f^0_1,\ldots ,f^0_k),(0:\cdots :0:1:0))$.
The image of $F$ is contained in the vector subspace $\mathbf{E}\subset\Mat_{n+1,k}(\C)\times\C^{n+1}$
defined by (\ref{matrix}), and it does not intersect
$\Y_{\dd,n}$, since we have assumed $(f_2,\ldots,f_k)\not\in V_{\dd',n-1,\C P^i}$
for any $(f_2,\ldots,f_k)$ such that $(f_1,f_2,\ldots,f_k)\in U$. Hence, due to proposition \ref{matrixlemma},
the intersection multiplicity of $F(U\times U_1)$ and $\W_{\dd,n}$ is equal to the intersection multiplicity of $F(U\times U_1)$ and $\X_{\dd,n}$.

Recall that in the the proof of proposition \ref{recndn} we have introduced the isomorphism $R:\mathbf{E}\to (\Mat_{n,k}\times\C^n)\times\C^2$; the composition $R\circ F$ takes
$$((f_1,\ldots,f_k),(0:\cdots :0:z_1:\cdots:z_{n-i}:1:z_{n-i+1}))$$ to $$(F'((f_1,\ldots,f_k),y),d_1z^{d_1-1}_{n-i+1},z_{n-i+1}),$$
where $y=(0,\ldots ,0,z_1,\ldots,z_{n-i},1)^T$, and $F'$ is the mapping $U\times \C^n\to\Mat_{n,k}(\C)\times\C^n$ given by
$$((f_1,\ldots,f_k),y)\mapsto (df_1|_{y},\ldots,df_k|_{y},y).$$ Since $U\pitchfork V_{\dd,n-1,\C P^{n-i}}$, the image of $F'$ intersects $\W_{k,n-1}$ transversally, which implies
the lemma.$\clubsuit$

%The following lemma can be proven in a similar way (cf. the computation of $F_{\dd,k}^*(\lk_{[\X_{k,n}],\mathbf{E}})$ and $F_{\dd,k}^*(\lk_{[\Y_{k,n}],\mathbf{E}})$
%in the proof of proposition \ref{recndn}).
\begin{lemma}\label{ggamma}
If $k>1$, we have $\gamma_i^{\dd,n}=d_1$.
\end{lemma}
{\bf Proof.} We repeat with minor modifications the proof of the previous lemma.
Let $U\subset\Pi_{\dd',n-1}$ be a small disc transversal to $V_{\dd',n-1,\C P^{n-i}}$ at a smooth point
$(f_2^0,\ldots,f_k^0)$, and let $U_1$ be the neighbourhood of $(0:\cdots :0:1:0)$ in $\C P^{n-i+1}\subset\C P^n$ introduced in the proof of lemma \ref{delta}.
Recall that $\gamma_i^{\dd,n}$ was defined in a slightly different way in the cases $k=n+1$ and $k<n+1$.

Consider first the case $k<n+1$. Choose $f_1^0$ so that $(\{f_1^0\}\times U)\cap V_{\dd,n-1,\C P^{n-i}}=\varnothing$ (where $\C P^{n-i}$ is
given by (\ref{cpn-i})). Assume that $\Sing (f_2,\ldots,f_k)\cap\C P^{n-i}=\{(0:\cdots:0:1)\}$.
Clearly, $\gamma_i^{\dd,n}$ is equal to the intersection multiplicity of $U'=S_1^{\dd,n}(U)$ and $V_{\dd,n,\C P^{n-i+1}}$
(the projective subspace $\C P^{n-i+1}$ being given by (\ref{cpn-i+1})). We proceed then as in the proof of lemma \ref{delta}, except that this time we use
%$$\tilde V_{\dd',n-1,\C P^{n-i}}=\{((f_2,\ldots,f_k),x)\in\Pi_{\dd',n-1}\times\C P^{n-i}\mid x\in\Sing(f_2,\ldots,f_k)\}$$
%instead of $V_{\dd,n-1,\C P^{n-i}}$
the isomorphism $T$ (and not $R$) introduced in the proof
of proposition \ref{recndn} instead of $R$ (cf. ibid, the computation of $F_{\dd,k}^*(\lk_{[\Y_{k,n}],\mathbf{E}})$).

If $k=n+1$, then $\gamma_i^{\dd,n}$ is equal to the intersection multiplicity of $\{0\}\times\C^{n+1}\subset\C^{n+1}\times\C^{n+1}$ and
the image of the mapping $U\times U_1\to\C^{n+1}\times\C^{n+1}$ that takes $((f_2,\ldots, f_k),x),x=(0:\cdots :0:z_1:\cdots:z_{n-i}:1:z_{n-i+1})$ to
$$((f_2(\tilde x),\ldots,f_k(\tilde x),z^{d_1}_{n-i+1})^T,\tilde x).$$ It can be easily seen that $\nu=d_1$, which
completes the proof of the lemma.$\clubsuit$

\subsection{Explicit formulae for $\m_i^{\dd,n}$}
Let us summarise the results of sections \ref{sec3} and \ref{sec4}.
\begin{lemma}\label{midnexpl}
The pullback $\bb_i^{\dd,n}\in H^{2i-1}(\GL_{n+1}(\C))$ of any class
$\aaa^{\dd,n}_i\in H^{2i-1}(\Pi_{\dd,n}\setminus\Sigma_{\dd,n}),i=1,\ldots,n+1,$ under an orbit map is equal to $\m_i^{\dd,n}\cc^{n+1}_i$ with
\begin{equation}\label{midnformulae}
\begin{array}{c}
\mbox{$\m_i^{\dd,n}=N(\dd,n)+(-1)^{n-k+1}$, if $i\geq n-k+2$, and}\\ 
\mbox{$\m_i^{\dd,n}=N(\dd,n)+(-1)^{i+1}N(\dd,n-i)$, if $i\leq n-k+1$.}
\end{array}
\end{equation}
Moreover, $\m_i^{\dd,n}>0$, unless $\dd=(2)$.
\end{lemma}

{\bf Proof.} Due to lemmas \ref{highnonempty}, \ref{prevnonempty} and proposition \ref{Ndn}, we have
$\bb_i^{\dd,n}=\m_i^{\dd,n}\cc^{n+1}_i,i=1,\ldots,n+1$, and $\m_i^{\dd,n}>0$ for $\dd\neq (2)$, so it remains only to prove (\ref{midnformulae}).
If $k=n+1$ or $k=1$ (and $i$ is arbitrary),
these formulae follow immediately from lemma \ref{highnonempty}, proposition \ref{krainie} and lemmas
\ref{delta} and \ref{ggamma}. If $\dd$ and $n$ are arbitrary, and $i=n+1$, then (\ref{midnformulae}) is given by lemma \ref{highnonempty}.

Let us now fix an $i$. Assume that $k\geq i$ (otherwise replace $\dd$ by a sufficiently long sequence of the form
$(2,\ldots,2,d_1,\ldots,d_k)$). Consider the following three arrays of numbers:

\begin{equation}\label{mass1}\begin{array}{c|cccc}
\vdots&\vdots&\vdots&\vdots&\vdots\\
i+l&(d_k-1)^{i+l+2}&N((d_{k-1},d_k),i+l+2)&\cdots\cdots&N(\dd,i+l+k)\\
\vdots&\vdots&\vdots&\vdots&\vdots\\
i-2&(d_k-1)^{i}&N((d_{k-1},d_k),i)&\cdots\cdots&N(\dd,i+k-2)\\
\vdots&\vdots&\vdots&\vdots&\vdots\\
0&(d_k-1)^2&N((d_{k-1},d_k),2)&\cdots\cdots&N(\dd,k)\\
-1&d_k-1&d_{k-1}d_k-1&\cdots\cdots&d_1\cdots d_k-1\\
\hline
&1&2&\cdots\cdots&k
\end{array}
\end{equation}
\begin{equation}\label{mass2}
\begin{array}{c|cccc}
\vdots&\vdots&\vdots&\vdots&\vdots\\
i+l&(-1)^{i+1}(d_k-1)^{l+2}&(-1)^{i+1}(N(d_{k-1},d_k),2+l)&\cdots\cdots&(-1)^{i+1}N(\dd,l)\\
\vdots&\vdots&\vdots&\vdots&\vdots\\
i&(-1)^{i+1}(d_k-1)^2&(-1)^{i+1}N((d_{k-1},d_k),2)&\cdots\cdots&(-1)^{i+1}N(\dd,k)\\
i-1&(-1)^{i+1}(d_k-1)&(-1)^{i+1}(d_{k-1}d_k-1)&\cdots\cdots&(-1)^{i+1}(d_1\cdots d_k-1)\\
i-2&(-1)^{i+1}&(-1)^{i+1}&\cdots\cdots&(-1)^{i+1}\\
\vdots&\vdots&\vdots&\vdots&\vdots\\
-1&1&1&\cdots\cdots&1\\
\hline
&1&2&\cdots\cdots&k
\end{array}
\end{equation}
\begin{equation}\label{mass3}\begin{array}{c|cccccc}
\vdots&\vdots&\vdots&\vdots&\vdots&\vdots&\vdots\\
i+l&\m_i^{(d_k),i+l+1}&\m_i^{(d_{k-1},d_k),i+l+2}&\cdots&\m_i^{(d_{k-i+1},\ldots,d_k),2i+l}&\cdots&\m_i^{\dd,k+i+l}\\
\vdots&\vdots&\vdots&\vdots&\vdots&\vdots&\vdots\\
i-2&\m_i^{(d_k),i-1}&\m_i^{(d_{k-1},d_k),i}&\cdots&\m_i^{(d_{k-i+1},\ldots,d_k),2i-2}&\cdots&\m_i^{\dd,k+i-2}\\
i-3& &\m_i^{(d_{k-1},d_k),i-1}&\cdots&\m_i^{(d_{k-i+1},\ldots,d_k),2i-3}&\cdots&\m_i^{\dd,k+i-3}\\
\vdots&&&\ddots&\vdots&&\vdots\\
-1&&&&\m_i^{(d_{k-i+1},\ldots,d_k),i-1}&\cdots&\m_i^{\dd,k-1}\\
\hline
&1&2&\cdots&i&\cdots&k\\
\end{array}
\end{equation}
Let us denote the items of the arrays (\ref{mass1}), (\ref{mass2}) and (\ref{mass3})
with coordinates $(a,b)$
by $x_1(a,b)$, $x_2(a,b)$ and $x_3(a,b)$ respectively. We have already seen that $x_3(a,b)=x_1(a,b)+x_2(a,b)$, when $b=-1$, or $a=1$ or
$a+b=i-1$. Due to proposition \ref{recndn} and lemma \ref{prevnonempty}, the items of the arrays satisfy the recursive relation
$x_j(a,b)=(d_{k-a+1}-1)x_j(a,b-1)+d_{k-a+1}x_j(a-1,b),j=1,2,3$. Hence, any $x_3(a,b)$ is equal to $x_1(a,b)+x_2(a,b)$.

This proves the formulae (\ref{midnformulae}) for $\m_i^{\dd,n}$ such that $k=\mbox{the length of $\dd$}$ is $\geq i$, and in fact,
for any $\m_i^{\underline{e},l}$, where $\underline{e}$ is a sequence of the form
$(d_{k-j},\ldots,d_k)$, and $l\geq\mathop{\mathrm{max}}(i-1,\mathop{\mathrm{length}}(e)-1)=\mathop{\mathrm{max}}(i-1,j)$. Hence, the assumption $k\geq i$ that we made does not restrict the generality. The lemma is proven.$\clubsuit$

\section{Proofs of the theorems}\label{sec5}
%In this section we suppose $\dd\neq (2)$.
\subsection{Proof of theorem \ref{main}}
\begin{Prop}\label{finitestab}
Suppose $\dd\neq (2)$. Then the stabiliser $G\subset\GL_{n+1}(\C)$ of an element of $(f_1,\ldots,f_k)\in\Pi_{\dd,n}\setminus\Sigma_{\dd,n}$ is finite.
\end{Prop}
{\bf Proof.} Denote by $H$ the connected component of the identity of $G$. To prove the proposition, it suffices to show that $H$
is trivial. If $k=n+1$, this follows from the fact that any element of $H$ acts identically on the preimage of a generic point $\in\C^{n+1}$ under
the ramified covering $\C^{n+1}\to\C^{n+1}$ defined by $x\mapsto (f_1(x),\ldots,f_{n+1}(x))^T$. 

Suppose now that $k<n+1$. Notice that the elements of $H$ diagonalise simultaneously in some basis of $\C^{n+1}$.
Indeed, if $(\dd,n)\neq ((3),2),((2,2),3)$, this follows easily from the absence of nonzero holomorphic vector fields on a smooth complete intersection
of multidegree $\dd$ in $\C P^n$ (see, e.g., \cite[proposition 2.11]{wahl}\footnote{I am grateful to J. Steenbrink for this reference.}); the remaining two cases
(which correspond to elliptic curves in $\C P^2$ or on a quadric in $\C P^3$) can be treated directly.

Hence, $H$ is in fact a complex torus. If $H\neq\{\mathrm{Id}\}$, the rational cohomology mapping induced by
$\GL_{n+1}(\C)\to\GL_{n+1}(\C)/H$ (and hence, the rational cohomology mapping induced by
$\GL_{n+1}(\C)\to\GL_{n+1}(\C)/G$) is not surjective, which contradicts lemma \ref{midnexpl}. $\clubsuit$

It should be not very difficult to show directly that the stabilisers of the elements of $\Pi_{\dd,n}\setminus\Sigma_{\dd,n}$ do not contain unipotent transformations, which would enable one to prove proposition \ref{finitestab} without using the nonexistence of holomorphic vector fields on smooth complete intersections.

The first part of theorem \ref{main} (the existence of the geometric quotient) follows from proposition \ref{finitestab}
and from the fact that $\Pi_{\dd,n}\setminus\Sigma_{\dd,n}$ is a hypersurface complement (and hence, an affine variety).

The second part
of theorem \ref{main} follows now from the Leray-Hirsch principle and lemma~\ref{midnexpl}.$\clubsuit$

%{\bf Remark.} The normal bundle of a complete intersection $X$ of multidegree $\dd$ in $\C P^n$ is isomorphic to (the restriction of)
%$$\bigoplus_{i=1}^k\xi^{d_i},$
%where $\xi$ the cotautological bundle. This implies that the Euler characteristic of $X$ is
%$$d_1\cdots, d_k\sum_{\begin{array{c}(\alpha_1,\ldots,\alpha_k)\\ \alpha_i\geq 0,\sum\alpha_i\leq n-1\end{array}}} (-1)^{\sum\alpha_i} C_{n+1}^{n-k-\sum\alpha_i}
%d_1^{\alpha_1}\cdots d_k^{\alpha_k}.$$

\subsection{Proofs of theorems \ref{upperboundvector} and \ref{upboundprojective}}
Both proofs are based on the following trivial observation.

\begin{Prop}\label{obs}
Suppose that $\rho:GL_l(\C)\to\GL_N(\C)$ is a representation that takes scalar matrices to scalar matrices.
Let $X\subset\C^N$ be a connected open subset invariant both under $\rho(\GL_l(\C))$ and $\C^*$.
Suppose that all stabilisers of both the action of $\GL_l(\C)$ on $X$ and the action of $\PGL_l(\C)$ on $X/\C^*$ are finite.
\begin{enumerate}\item Let $a_1,\ldots,a_l$ be cohomology classes of $X$ such that the pullback of any $a_i$ under an orbit map is $m_i\cc^l_i$ with $m_i\neq 0$. Then the order of the stabiliser of any
$x\in X$ divides $\prod_{i=1}^l m_i$.
\item If moreover there exist nonzero $u_2,\ldots,u_l$ such that any $u_i a_i$ descends into $X/\C^*$, then $l$ divides $\prod_{i=2}^l u_im_i$, and
the stabiliser of any $\bar x\in X/\C^*$ divides $\frac{1}{l}\prod_{i=2}^l u_im_i$.
\end{enumerate}
\end{Prop}
$\clubsuit$

The first assertion of this proposition together with lemma \ref{midnexpl} imply the following theorem
\begin{theoremx}
Let $G\subset\GL_{n+1}(\C)$ be the stabiliser of an element $(f_1,\ldots,f_k)\in \Pi_{\dd,n}\setminus\Sigma_{\dd,n}$. If $k=n+1$, then the order of $G$ divides
$$\prod_{i=1}^{n+1}(N(\dd,n)+1)=(d_1\cdots d_k)^{n+1};$$
if $k<n+1$, then the order of $G$ divides
$$\left(\prod_{i=1}^{n-k+1}(N(\dd,n)+(-1)^{i+1}N(\dd,n-i))\right)\left(\prod_{i=n-k+2}^{n+1}(N(\dd,n)+(-1)^{n-k+1})\right).$$
\end{theoremx}
$\clubsuit$

Recall that explicit expressions for $N(\dd,n)$ are given by proposition \ref{krainie} and lemma \ref{ladders}.

If $d\geq 3$, then, due proposition \ref{krainie}, we have $N((d),i)=(d-1)^{i+1}$. By substituting this in
theorem \thevb${}^{\prime}$, we obtain immediately the assertion of theorem \ref{upperboundvector}
(we make the change of variable $i\to n+1-i$ in the product to make it look nicer).$\clubsuit$
%Due to proposition \ref{obs}, order of the stabiliser of an $f\in\Pi_{\dd,n}\setminus\Sigma_{\dd,n}$ divides
%$$\prod_{i=1}^{n+1}\m_i{(d),n}.$$
%Our first task in the proof of both theorems \ref{upperboundvector} and \ref{upboundprojective} will be obtaining
%explicit expressions for $\bb^{\dd,n}_i,i=1,\ldots,n+1$ when $
%=(d)$.
%
%It follows from lemma \ref{midnexpl}
%\ref{highnonempty} and proposition \ref{krainie}
%that $\m^{(d),n}_{i}$ is equal to
%\begin{equation}\label{explhi}
%$$(N((d),n)+(-1)^n(d-1)^{n-i+1)=((d-1)^{n+1}+(-1)^n)\cc^{n+1}_{n+1}.$$
%\end{equation}
%(Actually, this formula for $\bb^{(d),n}_{n+1}$ is easy to obtain directly.)
%Lemma \ref{prenonvempty} allows us to express $\bb^{(d),n}_i,i\leq n$ as $\bb^{(d),n-1}_i\delta^{(d),n}_i$, but
%
%On the other hand, lemma \ref{prevnonempty}, proposition \ref{krainie} and lemma \ref{delta} imply that
%It follows from lemma \ref{midnexpl} that
%for any $i=1,\ldots,n+1$ we have
%$$\bb^{(d),n}_i=((-1)^{i+1}+(d-1)^{i})(d-1)^{n-i+1}\cc^{n+1}_i,$$ which implies immediately theorem \ref{upperboundvector}

Passing from the vector case to the projective one requires a little more work. It is only here that we make use of the general definition of linking numbers
given in \ref{BM1}.
%{\bf Proof of lemma \ref{explicit}.} Let us fix an $i\in\{ 1,\ldots,n\}$. Define the {\it iterated suspension map} $S:\Pi_{(d),i-1}\to\Pi_{(d),n}$ by
%$$f\mapsto f+\sum_{j=i}^n x_j^d.$$
%Let $P$ be the projectivisation of the vector subspace
%$$\{(x_0,\ldots,x_n)\mid x_j=0\mbox{ for $j<i-1$}\}$$ of $\C^{n+1}$; notice that $P\cong \C P^{n-i+1}$.
%
%The preimage of $V_{(d),n,P}$ under $S$ is $V_{(d),i-1,\{x_0\}}$, where $x'=(0:\cdots:0:1)$.
%
%Now set $f=x_0^d+\cdots+x_{i-1}^d$.
%If we identify $\GL_i(\C)$ with the image of the canonical inclusion $\GL_i(\C)\subset\GL_{n+1}(\C)$, the mapping $S$ becomes
%$\GL_i(\C)$-equivariant. 
%Due to proposition~\ref{lk},
%\begin{equation}\label{yyy}
%S^*(\lk_{V_{(d),n,P},\Pi_{(d),n}})=\mu\lk_{V_{(d),i-1,\{x'\}},\Pi_{(d),i-1}},
%\end{equation}
%where $\mu$ is the intersection multiplicity of $V_{(d),n,P}$ and the image of $S$.
%
%Notice that $V_{(d),i-1,\{x_0\}}$ is just a vector space;
%set $E$ to be the following vector subspace
%of $\Pi_{(d),i-1}$ transversal to $V_{(d),i-1,\{x'\}}$:
%$$\{\sum_{j=0}^{i-1}a_j x_{i-1}^{d-1}x_j\}.$$ The multiplicity $\mu$ is equal to the intersection multiplicity of $S(E)$ and $V_{(d),n,P}$ (the
%intersection consists of
%only one element $x_i^d+\cdots+x_n^d$).
%
%Lemma \ref{explicit} is obtained from (\ref{yyy}), (\ref{explhi}) and the following proposition:
%
%\begin{Prop}
%The intersection multiplicity of $S(E)$ and $V_{(d),n,P}$ is $(d-1)^{n-i+1}$.
%\end{Prop}
%

Let us introduce some additional notation. Let $m\leq n$ be an integer, and
set $$\x_{m,n}(d)=\frac{1}{m!}\left. \frac{d^m}{dt^m}\right|_{t=0}\frac{1}{(1+(d-1)t)^{n+1}}.$$ Let
$\y_{m,n}(d)$ be the right bottom item of the matrix $A^m$, where $$A={\begin{pmatrix}
0&0&\cdots&0&-\x_{m,n}(d)\\
1&0&\cdots&0&-\x_{m-1,n}(d)\\
\hdotsfor{5}\\
0&\cdots&1&0&-\x_{2,n}(d)\\
0&\cdots&0&1&-\x_{1,n}(d)
\end{pmatrix}}$$ Notice that, due to propositions \ref{chern} and \ref{degreeprop}, we have $\y_{m,n}(d)=\deg V_{(d),n,\C P^m}$.

\begin{Prop}\label{propdeg}
We have $\deg V_{(d),n,\C P^m}=\y_{m,n}(d)=C_{n+1}^m(d-1)^m$.
\end{Prop}
{\bf Proof.} Assume that $m\geq 2$ (otherwise, the assertion is obvious). Write $A^i$ as
$(a^i_{l,j})$. The proposition would follow, if we prove the following claim for any $i=1,\ldots,m$:

\medskip
{\bf Claim:} We have
$a^i_{m,m-i+l}=(d-1)^lC_{n+1}^l$ for $l=0,\ldots, i$, and $a^i_{m,j}=0$, if $j<m-i$.

\medskip
We proceed by induction on $i$.
The case $i=1$ is clear. Suppose that the claim is proven for some $i\leq m-1$. By writing $A^{i+1}=A^iA=A^i(N+A')$, where $N$ is nilpotent, and
all columns of $A'$ are zero,
except for the last one, which consists of $-\x_{m,n},\ldots,-\x_{1,n}$, one obtains
$$a^{i+1}_{m,j}=a^i_{m,j+1}=0 \mbox{ for $j<m-(i+1)$},$$
$$a^{i+1}_{m,m-(i+1)+l}=a_{m,m-i+l}^i=(d-1)^lC_{n+1}^l,l=0,\ldots,i,$$
$$a^{i+1}_{m,m}=-\sum_{j=m-i}^m a^i_{m,j}\x_{m-j+1,n}=-\sum_{l=0}^i (d-1)^lC_{n+1}^l\x_{i-l+1,n}.$$
Since $i+1\leq m$ and
$$(1+(d-1)t)^{n+1}(1+\sum_{j=1}^m\x_{j,n}t^j)=1+o(t^m),t\to 0,$$
we conclude that
$$\sum_{l=0}^i (d-1)^lC_{n+1}^l\x_{i-l+1,n}+(d-1)^{i+1}C_{n+1}^{i+1}=0,$$
which completes the proof of the claim.
$\clubsuit$

Now set $\Pi'_{(d),n}$ to be the projectivisation of $\Pi_{(d),n}$. Let $i$ be an integer such that $2\leq i\leq n+1$. Set
$\Sigma'_{(d),n}$ respectively, $V'_{(d),n,\C P^{n-i+1}}$, to be the image of $\Sigma_{(d),n}$,
respectively, of $V_{(d),n,\C P^{n-i+1}}$, in $\Pi'_{(d),n}$.

Let $a\in H^2(\Pi'_{(d),n})$ be the canonical generator, and set $b=[\Sigma'_{(d),n}]\frown a^{i-1}|_{\Sigma'_{(d),n}}$.
%Denote by $Q$ the intersection of $\Sigma'_{\dd,n}$ with a generic projective plane of dimension $\dim V'_{(d),n-i+1,\C P^{n-i+1}}+1$.
If $p=\deg V_{(d),n,\C P^{n-i+1}}=\deg V'_{(d),n,\C P^{n-i+1}},q=\deg\Sigma_{(d),n}=\deg\Sigma'_{(d),n}$, then the image of the class
$$c=\frac{1}{p}\LCM(p,q)[V'_{(d),n,\C P^{n-i+1}}]-\frac{1}{q}\LCM(p,q)b\in H_*(\Sigma'_{(d),n})$$
in $H_*(\Pi'_{(d),n})$ is zero, hence the linking number with this class can be viewed as a well-defined element of $H^*(\Pi'_{(d),n}\setminus\Sigma'_{(d),n})$.

Let $pr:\Pi_{(d),n}\setminus\{0\}\to\Pi'_{(d),n}$ be the natural projection.
Due to the first assertion of proposition \ref{lk}, the pullback of $\lk_{c,\Sigma'_{(d),n},\Pi'_{(d),n}}$ to $\Pi_{(d),n}\setminus\Sigma_{(d),n}$ is
the linking number with
$$\frac{1}{p}\LCM(p,q)[V_{(d),n,\C P^{n-i+1}}\setminus\{0\}]-\frac{1}{q}\LCM(p,q)pr^+_*(b)$$ in $\Pi_{(d),n}\setminus\{0\}$ (here $pr_*^+$ is the preimage map defined in \ref{BM1}). By proposition \ref{zzz} (in fact, this is the only time we use proposition \ref{zzz}),
$pr_*^+(b)=0$, and hence, the class $$\frac{\LCM(p,q)}{p}\lk_{[V_{(d),n}],\Pi_{(d),n}}=\frac{\LCM(C^{n-i+1}_{n+1}(d-1)^{n-i+1},(n+1)(d-1)^n)}{C^{n-i+1}_{n+1}(d-1)^{n-i+1}}
\aaa_i^{(d),n}$$ is
the pullback of an element of $H^*(\Pi'_{(d),n}\setminus\Sigma'_{(d),n})$. The proof of theorem \ref{upboundprojective} is now completed
using the second assertion of proposition \ref{obs} (as above, we make the change of variable $i\to n+1-i$ in the product).$\clubsuit$
%Denote by $\LCM(x,y)$ the least common multiple of the integers $x$ and $y$.
%Here is the table of $\x_{m,n}(d)$ and $\y_{m,n}(d)$ for $1\leq m\leq n\leq 4$:
%\begin{tabular}{|c|c|c|c|}
%\hline
%2(d-1)&&&\\
%3(d-1)&3(d-1)^2&&\\
%4*(d-1)&6(d-1)^2&4*(d-1)^3&\\
%
%\hline
%\end{tabular}

\section{Discussion}\label{discussion}
\subsection{Some particular cases of theorem \ref{upboundprojective}}

Here we list the values of (\ref{bound}) for $n=2,3,4$ and $3\leq d\leq 10$.

\begin{longtable}{|c|c|c|c|}
\caption{Values of (\ref{bound}) for small $n$ and $d$}\label{txxxx}\\
\hline
\backslashbox{$d$}{$n$}&2&3&4\endhead
\hline
3\rule{0pt}{11pt} & $432=2^4\cdot 3^3$ & $414720=2^{10}\cdot 3^4\cdot5$ & $218972160=2^{14}\cdot 3^5\cdot 5\cdot 11$ \\
\hline
4\rule{0pt}{11pt} &$18144=2^5\cdot 3^4\cdot 7$ & $2^{10}\cdot 3^8\cdot 5\cdot 7$&$2^{11}\cdot 3^{16}\cdot 5\cdot 7\cdot 61$\\
\hline
5\rule{0pt}{11pt} & $2^8\cdot 3\cdot 5^2\cdot 13$ & $2^{19}\cdot 3^2\cdot 5^3\cdot 13\cdot 17$ & $2^{30}\cdot 3^2\cdot 5^5\cdot 13\cdot 17\cdot 41$\\
\hline
6\rule{0pt}{11pt} & $2^4\cdot 3^3\cdot 5^4\cdot 7$ & $2^9\cdot 3^4\cdot 5^9\cdot 7\cdot 13$ & $2^9\cdot 3^5\cdot 5^{16}\cdot 7\cdot 13\cdot 521$\\
\hline
7\rule{0pt}{11pt} &$2^4\cdot 3^4\cdot 5\cdot 7^2\cdot 31$ & $2^{10}\cdot 3^8\cdot 5^2\cdot 7^3\cdot 31\cdot 37$ & $2^{14}\cdot 3^{16}\cdot 5^2\cdot 7^4\cdot 11\cdot 31\cdot 37\cdot 101$\\
\hline
8\rule{0pt}{11pt} &$2^7\cdot 3\cdot 7^4\cdot 43$ & $2^{13}\cdot 3^2\cdot 5^2\cdot 7^9\cdot 43$ & $2^{15}\cdot 3^2\cdot 5^2\cdot 7^{16}\cdot 11\cdot 43\cdot 191$\\
\hline
%\newpage
%\hline
9\rule{0pt}{11pt} &$2^{12}\cdot 3^5\cdot 7\cdot 19$ & $2^{28}\cdot 3^7\cdot 5\cdot 7^2\cdot 13\cdot 19$ & $2^{46}\cdot 3^9\cdot 5\cdot 7^2\cdot 11\cdot 13\cdot 19\cdot 331$\\
\hline
10\rule{0pt}{11pt} &$2^5\cdot 3^8\cdot 5^2\cdot 73$& $2^{11}\cdot 3^{17}\cdot 5^3\cdot 41\cdot 73$&$2^{11}\cdot 3^{32}\cdot 5^5\cdot 41\cdot 73\cdot 1181$\\
\hline
\end{longtable}

An explicit description of all automorphism groups of smooth projective hypersurfaces of given degree $>2$ is known in very few cases;
in fact, to the author's knowledge, there are three such cases: plane cubics and quartics and cubic surfaces.

A smooth cubic curve in $\C P^2$ can have 18, 36 or 54 projective automorphisms, depending on the value of the $j$-invariant. The least common multiple
of these numbers is $2^2\cdot 3^3$, which is 4 times smaller than the corresponding item of table \ref{txxxx}.

For $n=2,d=4$ the value of (\ref{bound}) is 18144, which is 9 times the LCM of the orders of the projective automorphism groups of smooth plane quartics
(the list of these groups is given e.g. in \cite[section 6.5.2]{dolg}\footnote{The reference was communicated to me by O. Tommasi.}).

The list of automorphism groups of smooth cubics in $\C P^3$ is given in \cite{hosoh} by T. Hosoh who
corrected an earlier classification by B. Segre \cite{segre}; the least common multiple of the orders of those groups
is $3240=2^3\cdot3^4\cdot 5$;
%$51870=2^7\cdot 3^4\cdot 5$ elements;
compare this with the $n=3,d=3$ item of table~\ref{txxxx}.
%which is $2^{10}\cdot 3^4\cdot 5$.

In fact, the expression (\ref{bound}) seems (at least for small $d$ and $n$) not to contain ``parasitic'' primes, i.e., primes that do not actually occur
as orders of automorphisms of smooth degree $d$ hypersurfaces in $\C P^n$.

\subsection{Odds and ends}
Here we consider some questions that do not fit into other sections, but are related to the topics discussed of the paper.
\subsubsection{Deck transformations}

The proof of the following statement is completely analogous to the proof of theorem \ref{upboundprojective}.
\begin{theorem}
Let $f:\C P^n\to\C P^n$ be a ramified covering of degree $dn$. Then the order of the group formed by the automorphisms $g:\C P^n\to\C P^n$ such that $f\circ g=f$ divides $$d^{n^2-1}\prod_{i=2}^{n+1}\frac{1}{C_{n+1}^i}\LCM(C_{n+1}^i,(n+1)d^{i-1}).$$
\end{theorem}
$\clubsuit$

This result should be easy to generalise to the case of ramified coverings from $\C P^n$ to arbitrary weighted projective spaces.
\subsubsection{Actions of the orthogonal groups}
Here we present two corollaries of theorem \ref{main} following from the well-known fact that for odd $l$ the inclusion $\SO_l(\R)\subset\GL_l(\C)$ induces an epimorphism of the rational cohomology groups (see, e.g., \cite{borel}).

%(We assume only that $q$ is not the square of a linear polynomial.)
%It is well known (see, e.g., \cite{borel}) that the inclusion $\SO_m(\R)\subset\GL_m(\C)$ induces an epimorphism of the rational cohomology groups
%for odd $m$. The following statement follows then easily from theorem \ref{main}.

\begin{corollary}
Let $q$ be a nonzero quadratic polynomial defining a (possibly singular) quadratic hypersurface $\subset\C P^n$.
Set
$m=\dim\ker q$. Suppose that $n+m$ is even and that $m\leq 1$.
The assertion of the division theorem \ref{main} remains true, if we replace $\GL_{n+1}(\C)$ by $\Aut^0 (q)$ and $\Pi_{\dd,n}\setminus\Sigma_{\dd,n}$ by
the space of all $(f_1,\ldots,f_k)$ of multidegree $\dd$ such that $\Sing(q,f_1,\ldots,f_k)=\varnothing$.
\end{corollary}
(Here $\Aut^0 (q)$ is the connected component of the identity of $\Aut (q)$; notice that $\Aut^0(q)$ contracts to $\SO_{n+1-m}(\R)\times\mathrm{U}_{m}$, hence the requirement
for $n+m$ to be even.)

$\clubsuit$

The condition that $n+\dim\ker q$ should be even can probably be removed.

Obvious as it is, this corollary indicates that one should be able to prove an analogue of theorem \ref{upboundprojective} for smooth hypersurfaces of quadrics. This will be the subject of a further work.

%Since for odd $m$ the mapping of the $H^*(\bullet,\Q)$ induced by the inclusion $\SO_m(\R)\subset\GL_m(\C)$ is epimorphic,
%we obtain the following statement.

\begin{corollary}
Denote by $\Pi_{\dd,n}(\R)\subset\Pi_{\dd,n}$ the set of fixed points of the (standard) complex conjugation, and suppose that $n$ is even.
The assertion of theorem \ref{main} remains true, if we replace $\Pi_{\dd,n}\setminus\Sigma_{\dd,n}$ by
$\Pi_{\dd,n}(\R)\setminus\Sigma_{\dd,n}$ and $\GL_{n+1}(\C)$ by
%the group
$\GL^+_{n+1}(\R)$.
\end{corollary}
$\clubsuit$

It is not known to the author if this statement holds for odd $n$.

\subsubsection{Possible generalisations and some open questions}
Notice that besides $\GL_{n+1}(\C)$ there are other groups acting naturally on $\Pi_{\dd,n}\setminus\Sigma_{\dd,n}$. For instance, set $\mathbf{G}_{\dd,n}$ to be the group generated by all transformations
$$(f_1,\ldots,f_i,\ldots,f_k)\mapsto (f_1,\ldots,af_i+gf_j,\ldots,f_k)$$
(where $i\neq j$ are indices such that $d_i\geq d_j$, $a\in\C^*$, and $g$ is a homogeneous polynomial of degree $d_i-d_j$).

The geometric quotient of $\Pi_{\dd,n}\setminus\Sigma_{\dd,n}$ by $\mathbf{G}_{\dd,n}$ obviously exists; its elements parametrise smooth complete intersections themselves, rather than their equations (since we quotient out all possible ways to pick up a minimal system of generators for
the homogeneous ideal of the variety given by $f_1(x)=\cdots=f_k(x)=0,(f_1,\ldots,f_k)\in\Pi_{\dd,n}\setminus\Sigma_{\dd,n}$).

Finally, notice that the subgroup of $\GL_{n+1}(\C)\times\mathbf{G}_{\dd,n}$ acting identically on $\Pi_{\dd,n}$ is isomorphic to $\C^*$. Denote by ${\cal G}_{\dd,n}$ the corresponding quotient group. The answers to the following questions are unknown to the author, if $k>1$.
\begin{itemize}
\item Do we have a division theorem for the action of $\mathbf{G}_{\dd,n}$ on $\Pi_{\dd,n}\setminus\Sigma_{\dd,n}$?
\item Suppose that $d_1=\cdots=d_k$. Do we have then a division theorem for the action of ${\cal G}_{\dd,n}$ on $\Pi_{\dd,n}\setminus\Sigma_{\dd,n}$? Notice that one can not expect such a result to hold for the sequences $\dd$ that contain at least two different items, since the cohomology of ${\cal G}_{\dd,n}$ will then have too many generators in dimension 1 (and $\Sigma_{\dd,n}$ is an irreducible hypersurface).
\item Notice that the group $\PGL_{n+1}(\C)$ acts on the quotient $(\Pi_{\dd,n}\setminus\Sigma_{\dd,n})/\mathbf{G}_{\dd,n}$; do we have a division theorem for that action?
\end{itemize}

\section*{Index of notation}
\nopagebreak

\begin{multicols}{2}
$k,$ \pageref{nkdd}

$n,$ \pageref{nkdd}

$\dd$, \pageref{nkdd}

$\Pi_{\dd,n},$ \pageref{PiSigma}

$\Sigma_{\dd,n}$, \pageref{PiSigma}

$\Sing(f_1,\ldots,f_k)$, \pageref{SING}

$\lk_{c,X,M},\lk_{c,M}$, \pageref{llk}

$\tot(\xi)$, \pageref{tot}

$\ee_m,$ \pageref{ecm}

$\cc_i^m,$ \pageref{ecm}

$o(\cc_i^m)$, \pageref{ecm}

$\Mat_{i,j}(\C)$, \pageref{mat}

$W_{k,n},$ \pageref{Wwkn}

$\W_{k,n}$, \pageref{Wkn}

$X_{k,n}$, \pageref{XY}

$Y_{k,n}$, \pageref{XXY}

$\X_{k,n}$, \pageref{XXYY}

$\Y_{k,n}$, \pageref{XXXYY}

$\aaa_i^{\dd,n},$ \pageref{a}

$\bb_i^{\dd,n},$ \pageref{b}

$\m_i^{\dd,n},$ \pageref{m}

$S_i^{\dd,n}$, \pageref{SSS}

$F_{\dd,n},$ \pageref{defndn}

$N(\dd,n)$, \pageref{defdefndn}

$\dd',$ \pageref{dd'}

$\gamma_i^{\dd,n}$, \pageref{deltagamma}

$\delta_i^{\dd,n}$, \pageref{deltagamma}
\end{multicols}

\begin{flushright}
{\sc Alexei Gorinov\\
Institute for Mathematics, Astrophysics and Particle Physics\\
Faculteit der Natuurwetenschappen, Wiskunde en Informatica\\
Radboud Universiteit Nijmegen\\
The Netherlands}\\
{\tt a.gorinov@math.ru.nl}
\end{flushright}

\thebibliography{99}
\bibitem{af} P. Aluffi, C. Faber, ``Linear orbits of smooth plane curves'',
J. Algebraic Geom., 2, 1993, no. 1, 155--184, {\tt arXiv:math.AG/9206001}.
\bibitem{borel} A. Borel, ``Sur la cohomologie des espaces fibr\'es principaux et des espaces homog\`enes de groupes de Lie compacts'',
Ann. of Math. (2) 57, 1953, 115--207.
\bibitem{brvet} W. Bruns, U. Vetter, ``Determinantal rings'', Lecture Notes in Mathematics, 1327, Springer-Verlag, Berlin, 1988.
\bibitem{dolg} I. V. Dolgachev, ``Topics in classical algebraic geometry, Part I'', available\footnote{As of
October the 24th, 2005.} at {\tt http://www.math.lsa.umich.edu/\symbol{126}idolga/topics1.pdf}
\bibitem{fuchs} A. T. Fomenko, D. B. Fuchs, ``A course in homotopic topology'', Nauka, Moscow, 1989 (in Russian).
\bibitem{quintics} A. G. Gorinov, ``Real cohomology groups of the
space of nonsingular curves of degree 5 in $\C P^2$'', Ann. Fac. Sci. Toulouse Math., 14, no. 3, 2005, 395-434, {\tt arXiv:math.AT/0105108}.
%\bibitem{gh}P. Griffiths, J. Harris, ``Principles of algebraic geometry'', Pure and Applied Mathematics, Wiley-Interscience (John Wiley \& Sons), New York, 1978.
\bibitem{hardy} G. H. Hardy, E. M. Wright, ``An introduction to the theory of numbers'', fourth edition, the Clarendon Press, Oxford University Press, London, 1975.
\bibitem{hosoh} T. Hosoh, ``Automorphism groups of cubic surfaces'',
J. Algebra 192, 1997, no. 2, 651--677.
\bibitem{hoso} A. Howard, A. J. Sommese, ``On the orders of the automorphism groups of certain projective manifolds'',
Manifolds and Lie groups (Notre Dame, Ind., 1980), 145--158, Progr. Math., 14, Birkh\"auser, Boston, Mass., 1981.
%\bibitem{heier} G. Heier, ``Effective finiteness theorems for maps between canonically polarized compact complex manifolds'', {\tt arXiv:math.AG/0311086}.
\bibitem{mm} H. Matsumura, P. Monsky, ``On the automorphisms of hypersurfaces'',
J. Math. Kyoto Univ., 3, 1963/1964, 347--361.
\bibitem{stepet} C.A.M. Peters, J.H.M. Steenbrink, ``Degeneration of the Leray spectral sequence for certain geometric quotients'',
Mosc. Math. J. 3, 2003, no. 3, 1085--1095, 1201, {\tt arXiv:math.AG/0112093}.
\bibitem{segre} B. Segre, ``The non-singular cubic surfaces,'' Oxford Univ. Press, London, 1942.
\bibitem{sz}E. Szab\'o, ``Bounding automorphism groups'', Math. Ann. 304, 1996, no. 4, 801--811.
\bibitem{tommasi} O. Tommasi, ``Rational cohomology of the moduli space of genus 4 curves'', Compos. Math. 141 200), no. 2, 359--384, 
{\tt arXiv:math.AG/0312055}.
%\bibitem{tsuji} H. Tsuji, ``Bound of automorphisms of projective varieties of general type'', {\tt arXiv:math.AG/0004138}.
\bibitem{vas1} V. A. Vassiliev, ``Geometric realization of the homology of classical Lie groups, and complexes that are $S$-dual to flag
manifolds'' (in Russian). Algebra i Analiz 3, 1991, no. 4,
113--120,
translation in St. Petersburg Math. J. 3, 1992, no. 4, 809--815.
\bibitem{vas2} V. A. Vassiliev, ``How to calculate homology groups of
spaces of nonsingular algebraic projective hypersurfaces in $\C P^n$'', Proc.
Steklov Math. Inst., 1999, 225, 121-140.
\bibitem{wahl} J. Wahl, ``Derivations, automorphisms and deformations of quasihomogeneous singularities'',
Singularities, Part 2 (Arcata, Calif., 1981), 613--624, Proc. Sympos. Pure Math., 40, Amer. Math. Soc., Providence, RI, 1983.
\bibitem{weis} B. Weisfeiler, ``Post-classification version of Jordan's theorem on finite linear groups''
Proc. Nat. Acad. Sci. U.S.A. 81, 1984, no. 16, Phys. Sci., 5278--5279.
\bibitem{xiao1}G. Xiao, ``Bound of automorphisms of surfaces of general type. I'',  Ann. of Math. (2) 139, 1994, no. 1, 51--77.
\bibitem{xiao2}G. Xiao, ``Bound of automorphisms of surfaces of general type. II'', J. Algebraic Geom. 4, 1995, no. 4, 701--793.
\end{document}